\documentstyle[draft]{article}
\title{Differential geometry of generalized almost quaterniohic structures, I}
\author{ V.F. Kirichenko, O.E. Arseneva}
\newcommand \implies{\Longrightarrow}
\newcommand \eq{\!\!&\!\!=\!\!&\!\!}
\newcommand \la{\langle}
\newcommand \ra{\rangle}
\newcommand \n[2]{\nabla_#1\,#2}
\newcommand \un[2]{({#1}_{#2})}
\newcommand \pq{\pi AQ_\alpha}
\newcommand \aq{\pi AQ_1}
\newcommand \ad{\mbox{\,\rm ad\,}}
\newcommand \Ad{\mbox{\,\rm Ad\,}}
\newcommand \tr{\mbox{\,\rm tr\,}}
\newcommand \End{\mbox{\rm End\,}}
\newcommand \id{\mbox{\,\rm id}}
\newcommand \re{\mbox{\rm Re\,}}
\newcommand \im{\mbox{\rm Im\,}}
\newcommand \procl[1]{\medskip\par{\bf#1.}\it}
\newcommand \eprocl{\medskip\par\rm}
\newcommand \definition[1]{\medskip\par{\bf#1.}}
\newcommand \edefinition{\medskip\par}
\newcommand \example[1]{\medskip\par{\bf#1.}}
\newcommand \eexample{\medskip\par}
\newcommand \remark[1]{\medskip{\it#1.}}
\newcommand \eremark{\par}
\newcommand \demo[1]{\par{\it#1.}}
\newcommand \edemo{\medskip\par}
\newcommand \qed{$\quad\Box$}
\newcommand \qeds{\quad\Box}
\renewcommand \a{\alpha}
\newcommand \e{\varepsilon}
\renewcommand \b{\beta}
\renewcommand \c{\gamma}
\renewcommand \d{\delta}
\renewcommand \l{\lambda}
\renewcommand \o{\overline}
\begin{document}
\maketitle
\begin{abstract} The fibre bundles adjoint to generalized almost quaternionic structures
are studied. The most important classes of generalized almost quater\-nionic
manifolds are considered.
\end{abstract}
\tableofcontents
\newpage
The important role of almost quaternionic structures in the row of
differen\-tial-geometric structures studied at present is explained primarily
by the fact that geometry of almost quaternionic manifolds generalizes the
most essential features of geometry of 4-dimensional oriented (pseudo-)Riemannian
manifolds being of great importance in theoretical and mathematical physics.
Among the features we shall stress the existence of self-dual and anti-self-dual
2-forms on 4-dimensional oriented (pseudo-)Riemannian manifolds generating
canonical almost quaternionic structures on the manifolds. Geometry of almost
quaternionic manifolds, in its turn, is closely related to Einsteinian geometry
of manifolds, the latter having been investigated by many outstanding geometers
[1]. Another important feature of 4-dimensional (pseudo-)Riemannian geometry
is connected with the remarkable discovery of twistor method by R.Penrose.
The method was proposed in the mid$-60s$ for solving problems in gravitation
theory and appeared to be fruitful in the field theory of Yang-Mills. The
famous Penrose's Twistor Programme [2] consists in using twistor correspondence
for transforming conformal invariant fields given on Minkowski complex space
subsets into objects of complex algebraic geometry (such as holomorphic bund\-les,
cohomologies with coefficients in analytic bundles and others) that are defined
on twistor space subsets. In the works by Penrose and other authors the programme
was used for gaugh fields, i.e. solving mass-free equations, on Minkowski spaces.
Interpretation of duality equations in terms of algebraic bundle over ${\bf C}P^3$
with the help of twistor transformations that was mentioned by Ward in 1977
allowed to reduce the problem of classification of instanton solutions to the
problem of algebraic geometry. In this way, for instance, there was received
the well-known Atiyah-Drinfeld-Hitchin-Manin classification (ADHM-construction)
[3], monopol classification of Nam-Hitchin [4], and others. Twistor construction
is closely connected with the existence of canonical almost quaternionic
structures on 4-dimensional oriented (pseudo-)Riemannian mani\-folds mentioned
above. Besides, S.Salamon [5] and L.B\'erard-Bergery [6] showed indepen\-dently that  the  construction  is
naturally transferred onto quaternionic manifolds of random dimension. It
allows to reduce the study of quaternionic manifolds to the study of complex
manifolds -- associated twistor bundle spaces. Using twistor construction
S.Salamon proved [7] that compact quaternionic-Kaehler manifold of positive
Ricci curvature for which the second Stiefel-Whitney class of bundle of purely
imaginary quaternions is equal to zero is isometric to canonical quaternionic
projective space ${\bf H}P^n$.
\par
The notion of almost quaternionic structure was first found in 1951 in P.Liber\-mann's
paper [8], the structure being considered in its narrow sense, as a pair of
anti-commuting almost complex structures on manifold. We call such structures
almost quaternionic with parallelizable structural bundle $(\pi AQ$-structures).
In the paper mentioned above as well as in [9] there was solved one of the most
important problems of differentional-geometric structure theory, viz. the problem
of integrability of $\pi AQ$-srtucture. Namely, P.Libermann showed that a
$\pi AQ$-manifold is integrable if and only if it is locally isomorphic to
quaternionic affine space. A great contribution to the study of
 $\pi AQ$-structures was made by M.Obata [10] who investigated the structure
of affine connections preserving $\pi AQ$-structure in parallel translations,
as well as properties of transformations  preserving $\pi AQ$-structure and by
S.Ishihara who studied special types of $\pi AQ$-manifolds transformations [11].
\par
Soon the reseachers realized that the class of $\pi AQ$-structures is too narrow
for constructing a self-contained  quaternionic geometry: even a quaternionic
projective space ${\bf H}P^n$ does not possess the structure of such type. So,
there was proposed (and it is generally accepted at present) a broader
understanding of almost quaternionic structure as a subbundle of tensor bundle
of type (1,1) on manifold whose type fibre is a quaternion algebra (or,
equivalently, as $GL(n,{\bf H}\cdot Sp(1)$-structures on manifold [12]). Some
authors called the structures almost quaternionic [13]. But the study of such
structures by traditional methods was difficult as the structures are not
strictly and globally defined by a given system of tensor fields like, for
example, Riemannian, or almost Hermitian, or almost contact structures. Thus,
the question of integrability of almost quaternionic structures is not answered
at present (unlike $\pi AQ$-structures). Moreover, it is not clear yet in what
terminology the integrability criterion could be formulated. Nevertheless,
several interesting results in this direction were received, for example,
Kulkarni theorem asserting that a compact simply-connected integrable quaternionic
manifold is isometric to the quaternionic projective space [14].
\par
The most important results in the theory of almost quaternionic manifolds were
received for quaternionic, quaternionic-Hermitian, quaternionic-Kaehler and
hyper-Kaehler manifolds. Quaternionic manifolds were first considered by
S.Salamon [15]. They are quaternionic counterpart of complex manifolds.
Nevertheless, their geometry differs considerably from complex geometry, i.e.
unlike Kaehler structures, a quaternionic-Kaehler structure is not always
inte\-grable. The basic property of quaternionic manifolds is integrability of a
canonical almost complex structure on the space of their twistor bundle for
manifolds of dimension
greater than 4 [7] (in case of 4-dimension its integrability is known to be
equivalent to manifold self-duality [16]). M.Berger proved [17] that
quaternionic-Kaehler manifold of dimension greater then 4 is an Einsteinian
manifold. It is Ricci-flat if and only if it is locally hyper-Kaehler. Otherwise,
it is not even locally reducible. The author also showed [18] that compact
oriented quaternionic-Kaehler manifold with positive sectional curvature is
isometric to the canonical quaternionic projective space. Hyper-Kaehler structures
were thoroughly investigated by A.Bea\-quville [19] who stated their close connection
with complex-symplectic structures. At present a collection of examples of
hyper-Kaehler manifolds of a random dimension are known [19,20]. On the other
hand, we know very few examples of complete non-hyper-Kaehler quaternionic-Kaehler
manifolds, all the manifolds being homoge\-neous. Moreover, D.V.Alekseevskii
proved [21] that any compact homogeneous quaternionic-Kaehler manifold is a
Riemannian symmetric space. Such spaces are completely classified by J.Wolf
[22]. The question of the existence of non-symmetric  quaternionic-Kaehler
manifolds was an open one. The first examples of the manifolds were given by
D.V.Alekseevskii [23].
\par
As mentioned above S.Salamon and L.Berard-Bergery have found independ\-ently a
new approach to the study of (almost) quaternionic manifolds genera\-lizing
Penrose twistor construction and allowing to reduce the study of quater\-nionic
manifolds to some problems of complex manifold theory. Within the approach
the authors received, for example, the following graceful and funda\-mental
characteristics of quaternionic-Kaehler manifolds:
\procl{Theorem} {\rm(Salamon [7,5], B\'erard-Bergery [6])} Let $M$
be a quaternionic-Kaehler manifold of positive Ricci curvature. Then its
twistor space $Z$ admits Kaehler-Einsteinian metric of positive Ricci curvature.
With respect to this metric the natural projection $\pi:Z\rightarrow M$ is a
Riemannian submersion with totally geodesic fibres. In particular, the compact
quaternionic-Kaehler manifold of positive Ricci curvature is simply-connected,
all its odd Betti numbers being equal to zero.
\eprocl
Unfortunately, with this approach the questions connected with possibility of
generalization on the multidimensional case of self-dual and anti-self-dual
2-form notions being basic in 4-dimensional Riemannian geometry have not been
discussed up to now. In particular, the possibility of generalization of
fundamental self-dual and anti-self-dual manifolds notions on the case of almost
quaternionic manifolds of arbitrary dimension has remained unclear. The role
of spinors in geometry of almost quaternionic manifolds has not been studied,
neither the interrelation of almost quaternionic structures in the sense of
Libermann (very convenient in applying traditional methods) and the accepted
broader understanding of almost quaternionic structures mentioned above.
\par
In the present paper the problems mentioned are investigated and as well as
the role of geometry of generalized almost quaternionic structures as generalization
of three-web geometry also widely discussed at present (see, for example,
[24,25]). Namely, the papers on three-web theory give many examples of three-webs
having a number of remarkable properties and showing close relations of the
theory with other areas of mathematics, such as, algebraic geometry [26],
quasi-group and loop theory [27], etc. This gives hope that the theory of
almost quaternionic structures of hyperbolic type, or almost antiquaternionic
structures, generalizing directly the three-web theory and being part of the
theory of generalized almost quaternionic structures studied here will become
the object of further serious investigation.
\par
Sections 1 and 2 introduce the notion of generalized almost quaternionic\linebreak
$(AQ_\alpha-)$ structure generalizing the notion of almost quaternionic structure
of classical type, as well as the notion of almost antiquaternionic structure,
in particular, the three-web structure, on the basis of generalized quaternion
algebra. We distinguish and investigate the notion of $AQ_\alpha$-structure of
spinor type and $AQ_\alpha$-structure with parallelizable structural bundle
and study the interrelation between the structure types. The spinor aspects of
geometry of spinor type $AQ_\alpha$-structures are also discussed, in particular,
spintensorial algebra of spinor type $AQ_\alpha$-manifold is constructed.
\par
In Section 3 we innerly add a special connection generalizing Chern con\-nection
of three-web theory to every $AQ_\alpha$-structure with parallelizable srtuct\-ural
bundle. In terms of curvature and torsion tensors of the connection we find
integrability criteria of $AQ_\alpha$-structure with parallelizable structural
bundles, as well as of its structural endomorphisms and basic distributions.
The notion of isoclinic distribution and isoclinic $AQ_\alpha$-structure is
introduced, and the criterion of isoclinic semi-holonomic $AQ_\alpha$-sructure
with parallelizable structural bundle, generalizing the well-known criterion
on three-web being isoclinic of M.A.Akivis, is found. We also introduce the
notion of isoclinic-geodesic $AQ_\alpha$-structure with parallelizable structural
bundle and find criteria of $AQ_\alpha$-structure with parallelizable structural
bundle being isoclinically geodesic.
\par
The developed theory allows us to construct whole classes of new interesting
(as we hope) examples of almost Hermitian and Einsteinian manifolds on a base
of canonical $\pi AQ_\alpha$-structure on a pull-back of certain manifolds.
Explicitly, we show that the pull-back of locally homogeneous naturally reductive
manifold carries two-parameter family of almost Hermitian structures of class
$G_1$ in Gray-Hervella classification [28]. The family contains the unique
quasi-Kaehlerian structure which is turn out to be nearly Kaehlerian structure.
Moreover, a corresponding two-parameter family of pseudo-Riemannian metrics
on such manifold contains exactly four Einsteinian metrics, if the manifold is
a semi\-simple Lie group equipped with Killing metric.
\par
We shall continue our research in the second part of this article. This section
will be devoted to applications of the developed methods to needs of self-dual
geometry and its multidimensional generalization. We shall introduce and
investigate the so called vertical type $AQ_\alpha$-structures which play
fundamental role in multidimen\-sional generalization of self-dual geometry. We
shall prove that quaternionic-Kaehler structures are the vertical type
$AQ_\alpha$-structures. We shall introduce and investigate generalized quaternionic
Kaehler manifolds and show that the manifolds of dimension greater than four
are Einsteinian, which generalizes the well-known result of M.Berger [17]. We
shall prove the multidimensional generalization of the well-known
Atiyah-Hitchin-Synger criterion of 4-dimensional oriented Riemannian manifold
being Einsteinian. We shall introduce the notions of twistor curvature tensor
and $t$-conformal semiflat $AQ_\alpha$-manifold and develop multidimensional
gene\-ralization of classical self-dual geometry on this base. Finally, we shall
consider 4-dimensional conformal semiflat generalized Kaehler manifold and
obtain some exhausting results in this direction.
\section{Spinor Geometry of Generalized Quaternion Algebra}
The characteristic feature of generalized almost quaternionic
structures is that in general case giving of the structure
on a manifold is not reduced to giving of one or several tensor
fields. This fact makes it difficult to directly apply classical
methods to its studying from the view point of differential geometry.
On the other hand, giving such a structure on a manifold
determines the fibre bundle over the manifold, whose standard
fibre being the algebra of generalized quaternions. Studying geometry
of generalized almost quaternionic structures is closely
connected with studying the orbits of linear representation of
such algebra. In particular, studying the natural representation
of the algebra being us to an interpretation of spinor geometry
on the basis of generalized quaternion algebra.
\subsection{Algebra of $\alpha$-quaternions}
Let {\it A} be an algebra over the field $F$, $char F\ne 2$, with
unit 1 and involution $a\rightarrow\o a \ (a\in A)$, and $a+\o a\in F,\ a\cdot\o a\in F$.
Recall [9] that the Cayley-Dickson duplication procedure allows
to construct by every nonzero element $\alpha\in F$ new algebra $(A,\alpha)$
which is $F$-modul $F^2=F\oplus F$, with the operation of multiplication
$(a_1,a_2)(b_1,b_2)=(a_1\,b_1+\alpha b_2\,\o a_2,\o a_1\,b_2+b_1\,a_2)$.
Here {\it A} is enclosed in $(A,\alpha)$ as a subalgebra of pairs of the type $\{\,(a,0)\mid a\in A\,\}$.
Denote $i=(0,1)\in (A,\alpha)$, then $i^2=\alpha$. Now, by the above
identification we have $(A,\alpha)=A\oplus iA$. Note that $(A,\alpha)$ is also
an algebra over $F$ with involution: if $x=a_1+ia_2\in (A,\alpha)$,
then $\o x=\o a_1-ia_2$. The involution has the same properties: $x+\o x\in F,\ x\cdot\o x\in F
\quad (x\in (A,\alpha))$. Moreover, if the quadratic form
$|a|^2=a\cdot\o a$ is non-degenerate on $A$, then the quadratic form
$|x|^2=x\cdot\o x$ is also non-degenerate on $(A,\alpha)$\ (see [9], p.42).
This involution is called {\it the operation of conjugation.}
\par
If $A=F=\bf R$, then $({\bf R},\alpha)={\bf R}\oplus i{\bf R}$ up to multiplication of
$i$ by constant multiplier $\sqrt{|\,\alpha|}$ is either the field $\bf C$ of complex
numbers $(\alpha =-1)$ or the ring $\bf D$ of double numbers $(\alpha =1)$.
Denote $({\bf R},\alpha)$ by ${\bf K}_\alpha$.
\par
Applying Cayley-Dickson procedure to ${\bf K}_\alpha$ we get the algebra
$(({\bf R},\alpha),\beta)=\bf R\oplus\,i{\bf R}\oplus\,j{\bf R}\oplus\,k\bf R$,
(where $j=(0,1)\in ({\bf K}_\alpha,\beta),k=i\cdot j)$, called {\it the algebra
of generalized quaternions} which is known to be associative and non-commutative [9, p.44]. This
algebra up to multiplication of $i$ and $j$ by constants $\sqrt{|\,\alpha|}$ and $\sqrt{|\beta|}$
(called {it calibration}), respectively, is either the body $\bf H$
of quaternions $(\alpha =\beta =-1)$ or the ring $\bf H\,^\prime$ of antiquaternions
$(\alpha =\beta =1)\ ($the case $\alpha\cdot\beta <0$ corresponds to the algebra of
antiquaternions). As the calibration does not make any essential
changes in the properties of generalized quaternion algebra we
consider the algebra calibrated and denote ${\bf H}_{\alpha}=(({\bf R},\alpha ),\alpha); \alpha=\pm 1,$
We call it {\it the algebra of} $\alpha $-{\it quaternions.} Note that for
$\alpha$ -quaternion $q=a+ib+jc+kd=(a+ib)+j(c-id)$ we have $k^2=-1;
\o q=a-ib-jc-kd; \ |q|^{2}=q\cdot\o q=a^{2}-\alpha b^{2}-\alpha c^{2}+d^{2}$.
\par
It follows that any $\alpha$-quaternion $q \in \,{\bf H}_\alpha$ admits a representation
of the type $q=z_1+jz_2;\quad z_1, z_2\in \,{\bf K}_\alpha,\ j^2=\alpha$ allowing
to construct in the natural way a representation of algebra ${\bf H}_\alpha$
into algebra $\End{\cal S}$ of endomorphisms of algebra ${\bf H}_\alpha$ regarded as
${\bf K}_\alpha$-modul ${\cal S}={\bf K}_\alpha \,\oplus{\bf K}_\alpha$ that we call
{\it an} $\alpha${\it-spinvector space.}
Namely, let $q=z_1+jz_2\in{\bf H}_\alpha, X=(X_1,X_2)\in{\cal S}$.
Put $[q]X=q\cdot (X_1+jX_2)=(z_1\,X_1+\alpha \o z_2\,X_2)+j(z_2\,X_1+\o z_1\,X_2)$.
Evidently, the matrix of endomorphism $[q]$ in basis $\{1,j\}$ has form:
\begin{equation}
(q)=\left(\begin{array}{rr}
z_1&\alpha\o z_2\\z_2&\o z_1\end{array}\right).
\end{equation}
\par
Note that the $\alpha$-spinvector space $\cal S$ is innerly endowed with
involution $\tau: X \rightarrow \o{X}$ generated by the operator of conjugation of
$\alpha $-quaternions and called {\it the operation of conjugation of} $\alpha ${\it -spinvectors.}
Besides, $\cal S$ is innerly endowed with  Hermitian  metric $\la\la\cdot,\cdot\ra\ra$.
Namely,  let $X=(z_1,z_2)\in{\cal S}$. Denote $z_1=\re\,X,\ z_2=\im\,X$.
Since $\o{X}=(\o{z_1},-z_2),\  \re \o X=\o{\re X},\ \im\,\o X=-\im X$. Further,
let $Y =(u_{1},u_{2})$ be another $\alpha $-spinvector. Put $\la\la X,Y\ra\ra
=\re(\o XY)=\o z_1\,u_1-\alpha\o z_2\,u_2$. Evidently, this is a non-degenerate
Hermitian form, $\la\la qX,qY \ra\ra=\re\,(\o{qX}\cdot qY)=\re\,(\o{X}\cdot
\o{q}\cdot q\cdot Y)=|q|^2\la\la X,Y \ra\ra$. Thus, we get
\procl{Proposition 1} Natural representation of non-isotropic really
normalized $\alpha\!\!~$-quaternions in the space $\cal S$ is realized
by conformal endomorphisms of module $\cal S$.\qed
\eprocl
Note that $|q|^2=q\cdot\o q=(z_1+jz_2)(\o z_1-jz_2)=z_1\,\o z_1
-\alpha z_2\o z_2=\det(q)=\det[q]$. In particular, unit $\alpha$-quaternions, that in
view of algebra ${\bf H}_\alpha$ being compositional form a multiplicative
subgroup of ring $\End{\cal S}$, are realized in this representation as
unitary unimodular endomorphisms of module $\cal S$.  Thus, we get
\procl{Proposition 2} Multiplicative subgroup of unit $\alpha$-quaternions
in natural repre\-sentation is identified with Lie group $SU(2,{\bf K}_\alpha)
=SP(1)$.\qed\eprocl
\remark{Remark} If $\alpha =1$ Lie group $SU(2,\bf D\,)$ is naturally
isomorphic to  Lie group $SL(2,{\bf R})$.  This isomorphism  allowing to
identify the given Lie groups juxtaposes the matrix
$$
\left(\begin{array}{cc}x_1+iy_1&x_2-iy_2\\x_2+iy_2&x_1-iy_1\end{array}\right)
\in SU(2,{\bf D}\,)
$$
to the matrix
$$
\left(\begin{array}{cc} x_1+y_1&x_2-y_2\\x_2+y_2&x_1-y_1\end{array}\right)
\in SL(2,\bf R).
$$
In particular, Lie algebras ${\sf su}(2,{\bf K}_\alpha)$ are semisimple.
\eremark
Lie group $SU(2,{\bf K}_\alpha)$ is a 3-dimensional Lie group acting on
its Lie algebra ${\sf su}(2,{\bf K}_\alpha)$ by a adjoined
representation $\Ad(g)X=gXg^{-1}\quad (g\in SU(2,{\bf K}_\alpha)$,
$X\in{\sf su}(2,{\bf K}-\alpha))$, orthogonal in the Killing
form of this Lie algebra. In view of the above and in force of
the connectivity of Lie group $SU(2,{\bf K}_\alpha)$ the image of mapping $\Ad$
lies in Lie group $SO(2-\alpha,1+\alpha;{\bf R})$. Thus, the homomorphism $s=\Ad:
SU(2,{\bf K}_\alpha )\rightarrow SO(2-\alpha ,1+\alpha ;{\bf R})$ is innerly defined.
In view of Lie group $SU(2,{\bf K}_\alpha )$ being semisimple the mapping
$s_{*}$ being a adjointed representation of its Lie algebra is non-degenerate and since
$\dim\,{\sf su}(2,{\bf K}_\alpha)=\dim\,{\sf so}(2-\alpha,1+\alpha;{\bf R})=3$,
it is an isomorphism of Lie algebra, and therefore, the mapping $s$-(a two-leaf one) is  a
covering mapping. The (two-valued) mapping $s^{-1}: SO(3,{\bf R}) \rightarrow  SU(2,{\bf C}\,)$
is a classicall spinor representation of an orthogonal group when $\alpha =-1$.
In the case of arbitrary algebra ${\bf H}_\alpha$ we follow the same notation.
\par Note that subspace $IQ\subset{\bf H}_\alpha$ of purely imaginary
$\alpha$-quaternions is cano\-nically identified with the space of Lie algebra
${\sf su}(2-\alpha,1+\alpha;{\bf R})$ that, in its turn,
admits canonical identification with Lie algebra
${\sf so}(2-\alpha,1+\alpha ;{\bf R})$
by means of mapping $s_{*}=\ad:{\sf su}(2,{\bf K}_\alpha)\rightarrow
{\sf so}(2-\alpha,1+\alpha;{\bf R})$. With the above identifications and in view of
(1) the purely imaginary $\alpha$-quaternion $q=ia+jb+kc$ has the
corresponding matrix
$$(q)=\left(\begin{array}{rr}ia&\alpha (b+ic)\\b-ic&-ia\end{array}\right)$$
which is an element of Lie  algebra ${\sf su}(2,{\bf K}_\a)$.
Also, $[(q_1),(q_2)]=(q_1)(q_2)-(q_2)(q_1)$. Now, in force
of the above identification  the space $IQ$  acquires the structure of Lie
algebra isomorphic to ${\sf su}(2,{\bf K}_\a)$ with commutator
$[q_1,q_2]=q_1\,q_2-q_2\,q_1$. Besides, $s_{*}(q)X=\ad(q)X=(q)X-X(q)$, hence we have:
$$s_{*}(q)=\left(\array{rrr}
0&2\alpha c&-2\alpha b\\-2\alpha c&0&2\alpha a\\-2b&2a&0\endarray\right)$$
in the basis $\{\,(i),(j),(k)\,\}$. In view of the above it is easy to
compute that $$-\tr(s_{*}(q)\,s_{*}(q))=-4\tr((q)(q))=8q\o q,$$
that is we get the valid
\procl{Proposition 3} The norm of purely imaginary $\alpha$ -quaternions
as elements of Lie algebras ${\sf su}(2,{\bf K}_\a)$ or
${\sf so}(2-\a ,1+\a;{\bf R})$ in Killing metric differs from their norms as
elements of algebra ${\bf H}_\a$ respectively only by the multiplier 2 or 8.\qed
\eprocl
\subsection{Spinbases}
The norm $|q|^2=q\cdot\o{q}$ of $\alpha $-quaternions allows to regard the
algebra ${\bf H}_\alpha$ as a real Euclidean space in which scalar product is
reconstructed by polarization:
\begin{equation}
\la q_{1},q_{2}\ra=\frac12(\o q_1\cdot q_2+\o q_2\cdot q_1);\qquad q_1,
q_2\in{\bf H}_\alpha.\end{equation}
Note that $\im(\o q_1\cdot q_2)=\im\o{(\o q_2\cdot q_1)}=-\im(\o q_2\cdot q_1)$,
hence $\la q_1,q_2\ra=\frac12(\re(\o q_1\cdot q_2)+\re(\o q_2\cdot q_1))$.
Recalling definition of the form $\la\la\cdot,\cdot\ra\ra$ we get
\begin{equation}
\la X,Y\ra =\frac12\{\,\la\la X,Y\ra\ra+\la\la Y,X\ra\ra\,\};\qquad
X,Y\in{\cal S}.\end{equation}
Further, it is evident that
\begin{equation}
\la\la X,Y\ra\ra =\la X,Y\ra +i\la X,i^3Y\ra;\qquad X,Y\in{\cal S}.
\end{equation}
Consider (2) again. In particular, let $q_1,q_2 \in IQ$. Then $\o q_1=-q_1,
\o q_{2}=-q_{2}$, and (2) assumes the form
\begin{equation}
q_1\cdot q_2+q_2\cdot q_1=-2\la q_1,q_2\ra;\qquad q_1,q_2\in IQ.
\end{equation}
Now we get the following
\procl{Proposition 4} The system $\{j_{1},j_{2},j_{3}\}\in IQ$ forms
an orthonormal basis of the space IQ iff
$$j_k\cdot j_m+j_m\cdot j_k=\left\{\begin{array}{rl}\pm 2,&\mbox{for $k=m$}\\
0,&\mbox{for $k\ne m.$}\end{array}\right.\qeds$$
\eprocl
For example, the system $\{i,j,k\}$  forms an orthonormal basis in $IQ$.
Recall that matrices having these elements as endomorphisms
of ${\bf K}_\alpha$-module {\cal S} have (by (1)) the following form:
\begin{equation}
\sigma_1=(i)=\left(\begin{array}{rr}i&0\\0&-i\end{array}\right);\
\sigma_2=(j)=\left(\begin{array}{rr}0&\a\\1&0\end{array}\right);\
\sigma_3=(k)=\left(\begin{array}{rr}0&\a i\\-i&0\end{array}\right);
\end{equation}
In is easy to see that when $\alpha=-1$ the matrices differ from the
classical Pauli ones only by multiplier $i$. We call them
{it generalized Pauli matrices}.
\procl{Theorem 1} For any orthonormal basis $\{\,j_{1},j_{2},j_{3}\,\}$
of the space IQ there exists an orthonormal basis $\{\,\e_1,\e_2\,\}$
of ${\bf K}_\a$-module{\cal S}, where $(j_{1})=\sigma_{1},(j_{2})=
\sigma_{2},(j_{3})=\pm\sigma_{3}$. The basis is  defined up to the action
of group $U(1,{\bf K}_\a)$ enclosed in $U(2,{\bf K}_\a)$ by diagonal matrices.
\eprocl
\demo{Proof} Consider vectors $\e'=\frac{X+i[j_1]^3X}{|X+i[j_1]^3X|}$
and $\e'=\frac{X-i[j_1]^3X}{|X-i[j_1]^3X|}$, where $X\in{\cal S}$ is
any vector not being eigenvector of operator $[j_1],\ i$ is imaginary unit
of ring ${\bf K}_\a$. Evidently, the vectors are eigenvectors of operator
$[j_1]$ having eigenvalues
$i$ and $-i$, respectively. In view of eigenspaces of the operator being
one-dimensional, each of the vectors is defined up to multiplication by
element $z\in{\bf K}_\a, |z|=1$. It is clear that in basis $\,\{\e'_1,
\e'_2\,\}\quad (j_{1})=\sigma _{1}$. Further, $[j_{1}]\e'_1=i\e'_1$ and hence
$[j_{2}][j_{1}]\e'_1= i[j_{2}]\e'_1$. On the other hand
$[j_{2}][j_{1}]\e'_1=-[j_{1}][j_{2}]\e'_1$, hence $[j_{1}][j_{2}]\e'_1
=-i[j_{2}]\e'_1$, i.e. $[j_{2}]\e'_1=a\e'_1,\quad a\in{\sf K}_\a.$
Similarly, $[j_{2}]\e'_2= b\e'_1,\quad b\in{\bf K}_\a$. Thus,
$$
(j_2)=\left(\begin{array}{cc}0&b\\a&0\end{array}\right),
$$
and since $(j_2)\in{\sf su}(2,{\bf K}_\a),\ \o a=ab$. Besides, since
$|j_2|^2=-\a$, i.e. $\det[j_{2}]=-\a$, we have $ab=\a,\quad
|a|^2=a\o a=\a ab=a^2=1,\ |b|^2=b\o b=\a\o a\o b=1$. Therefore,
$$
(j_2)=\left(\begin{array}{rr}0&\alpha\o a\\a&0\end{array}\right),\qquad |a|^2=1.
$$
Assume $\e_1=x\e'_1, \e_2=y\e'_2$,\ so that in basis $\{\e_1,\e_2\}$ we have
$(j_2)=\sigma_2$. Now we have: $[j_2]\e_1=x[j_2]\e'_1=xa\e'_2=xay^{-1}\e_2;
[j_2]\e_2=y[j_2]\e'_2=\a y\o a\e'_1=\a y\o ax^{-1}\e_1$.
Thus, $x$ and $y$  must be chosen so that $xy^{-1}=\a a^{-1}$.
Then, automatically, $\alpha y\o ax^{-1}=(xy^{-1})^{-1}\alpha\o a
=\alpha a\alpha\o a=\alpha^2|a|^2=1,\ xay^-1=xy^-1a=\a$. Thus, we can choose
an arbitrary element of ring ${\bf K}_\alpha$ with absolute value 1 as $x$,
and assume $y=\alpha ax$. Then in basis $\{\,\e_1,\e_2\,\}$
defined up to multiplication of its elements by $z\in{\bf K}_\alpha,\ |z|=1,
(j_2)=\sigma_2$. Finally, we similarly get that in this basis
$$
(j_{3})=\left(\begin{array}{rr}0&\alpha\o b\\b&0\end{array}\right),
\qquad b\in{\bf K}_\a,|b|^2=-\a.
$$
and in view of $\la j_2,j_3\ra =0$, and hence $[j_2][j_3]+[j_3][j_2]=0$, we have:
$$
\left(\begin{array}{cc}0&\alpha\\1&0\end{array}\right)
\left(\begin{array}{rr}0&\alpha\o b\\b&0\end{array}\right)\ +\
\left(\begin{array}{rr}0&\alpha\o b\\b&0\end{array}\right)
\left(\begin{array}{cc}0&\alpha\\1&0\end{array}\right)=0,
$$
hence $\alpha (b+\o b)=0,$ i.e. $b=\pm i$, and $(j_3)=\pm\sigma_3.$\qed
\edemo
\definition{Definition 1} An orthonormal basis of ${\bf K}_\a$-module {\cal S}
corresponding to the orthonormal basis of real space of purely imaginary
$\alpha$-quaternions is called {\it a spinbasis.}
\edefinition
In view of linear connection of Lie groups $SU(2,{\bf K}_\alpha)$ we get
\procl {Proposition 5} Orthonormal basis $\{\,j_1,j_2,j_3\,\}$ of the space
$IQ$ has the same orientation as the basis $\{\,i,j,k\,\}$ of the same space
iff in the spinbasis corresponding to the basis $\{\,j_1,j_2,j_3\,\}$ we have
$$(j_k)=\sigma_k;\qquad k=1,2,3.$$
\eprocl
\demo{Proof} Let $g\in SU(2,{\bf K}_\alpha)$ be the element transferring
spinbasis $\{\,\e^0_1,\e^0_2\,\}$, where $(i)=\sigma_1,(j)=\sigma_2,
(k)=\sigma_3$ into the spinbasis, where $(j_k)=\sigma _k, k=1,2,3;\omega(t)$
is the path into $SU(2,{\bf K}_\a)$ connecting unit $e\in SU(2,{\bf K}_\a)$
with $g$. Then $s\circ\omega(t)$ is the path in the group space of Lie group
$SO(2-\a,1+\a;{\bf R})$ deforming basis $\{\,i,j,k\,\}$ into
basis $\{\,j_1,j_2,j_3\,\}$.\qed\edemo
\subsection{Linear representation orbits of $\a$-quaternion algebra}
Let the precise linear representation of algebra ${\bf H}_\alpha$ be fixed.
Then the algebra is realized as a subalgebra of algebra $\End V$
of endomorphisms of (real) linear space $V$. Consider the orbit $Orb(X)=
\{\,Y\in V\mid Y=q(X);\ q\in{\bf H}_\alpha\}$ of the nonzero element $X\in V$.
It can be regarded as the image of homomorphism $X:{\bf H}_\a\rightarrow V;\
X(q)=q(X)$. Evidently,\ \ $\ker X\subset\{q\in{\bf H}_\alpha\mid |q|=0\,\}
\quad (X\in V)$.
This means that when $\ker X\ne \{0\}, X\in\ker q$, where $q$ is an isotropic
$\alpha$-quaternion. Thus, when $X\in V$ does not belong to the kernel of
isotropic $\alpha$-quaternion, $\ker X=\{0\}$ and $\dim Orb(X)=4$.
\par
In order to study the case when $X$  belongs to the kernel of isotropic
$\a$-quaternion we build the following construction. Consider ${\bf K}_\a$-module
$V^{{\bf K}_\alpha}=V\otimes_{\bf R}{\bf K}_\alpha$ and the naturally induced
linear representation of algebra ${\bf H}_\alpha$ on $V^{{\bf K}_\alpha}$
that we call {\it an extended representation} of the algebra. Let
$\{\id,I,J,K\}$ is orthonormal basis of space ${\bf H}_\alpha$. Evidently,
$V^{{\bf K}_\alpha}=D^i_I\oplus D^{-i}_I$,
where $D^{\lambda}_{I}$ is the eigensubmodule of operator $I$ with eigenvalue
$\lambda,i$ is the imaginary unit of ring ${\bf K}_\alpha; i^2=\alpha$. We have
$\sigma=\frac12(\id+iI^3)$ is a projector on $D^i_I$  and $J|D^i_I:D^i_I
\rightarrow D^{-i}_I$ is an isomorphism of ${\bf K}_\alpha$-modules. Let
$X\in D^i_I,\ \ q=a\id+bI+cJ+dK\in{\bf H}_\a$. Then $q(X)=(a+bi)X+(c-di)JX$;
$q(JX)=\alpha (c+di)X+(a-bi)JX$. Thus $Orb(X)={\cal L}(X,JX)$. In particular,
$\dim Orb(X)=2$, and
$$(q|Orb\,X)=\left(\begin{array}{cc}a+bi&\alpha(c+di)\\c-di&a-bi\end{array}
\right).$$
We get
\procl{Proposition 6} Orbits of eigenvectors of $\alpha$-quaternions
with purely imagi\-nary eigenvalues are 2-dimensional and have the natural
structure of $\alpha$-spin\-vector space.
\eprocl
Now, let $X\in V^{{\bf K}_\alpha}$ be an arbitrary element. Then $X=a+u,
a\in D^i_I,u\in D^{-i}_I$. Thus $Orb(X)={\cal L}(a,Ja,u,Ju)={\cal L}(a,Ju)
\oplus{\cal L}(Ja,u)$.
Indeed, ${\cal L}(a,Ju)\subset D^i_I,{\cal L}(Ja,u)\subset D^{-i}_I$, and since
$D^i_I\cap D^{-i}_I=\{0\}$, we have: ${\cal L}(a,Ju)\cap{\cal L}(Ja,u)=\{0\}$.
Thus, either $\dim Orb(X)=2$, or $\dim Orb(X)=4$, and $\dim Orb(X)=2$ iff
vectors $\{a,Ju\}$  are linearly dependent. Therefore, in the former case
$Orb(X)$ is the orbit of eigenvector $a$ (or $u$) having, as shown above,
a natural structure of $\alpha$ -spinvector space, and in the latter case
$Orb(X)$ is the direct sum of two such subspaces conjugated with respect to $J$.
So, we can have
\definition{Definition 2} 2-dimensional orbits of an extended linear representation
of $\alpha$-quaternion algebra are call {\it spinor subspaces}, its 4-dimensional
orbits being {\it bispinor subspaces} of space $V^{{\bf K}_\alpha}$.
\edefinition
\par
Consider, finally, the case when $X\in V$. Then $X=a+\tau a,\quad a=\sigma X$,
where $\tau:V^{{\bf K}_\a}\rightarrow V^{{\bf K}_\a}\quad$ is the operator
of natural conjugation. Thus, $Orb(X)={\cal L}(a,Ja,\tau a,J\tau a)=
{\cal L}(a,J\tau a)\oplus{\cal L}(Ja,\tau a)$. Here $\dim Orb(X)=2\quad$ iff vectors
$\{Ja,\tau a\}$ are linearly dependent. Let $Ja=\l\tau a,\quad\l\in{\bf K}_\a$.
We have $J^2a=\l J\tau a=\l\tau(\l\tau a)=\l\o{\l}a,\quad\l\in{\bf K}_\a$.
On the other hand, $J^2a=\alpha a$. Thus, $|\lambda|^2=\alpha$.
If ${\bf K}_\alpha={\bf C}$, then $\alpha=-1$ and this case is impossible.
Hence, ${\bf K}_\alpha={\bf D},\alpha=1,|\l|=1$. Besides, $\tau JA=\o{\l}a$.
But $\tau Ja=\frac12\tau(JX+iJIX)=\frac12(JX-iJIX)=\frac12(JX+iIJX)$, so
$JX+iIJX=\o{\lambda}(X+iIX)$, hence $JX=\o{\lambda}X$ and thus, $\lambda=\pm1$,
that is $X\in D^{\pm1}_{I}$. Inversely, if $X\in D^{\pm1}_{I}$, then
$q(X)=(a\pm b)X+(c\mp d)JX, q(JX)=(c\pm d)X+(a\mp b)JX$. In particularly,
$$
\dim Orb(X)=2,\quad (q)=\left(\begin{array}{cc}a\pm b&c\pm d\\c\mp d&a\mp b
\end{array}\right).
$$
Thus, the orbit of element $X\in V$ is 2-dimensional iff $X$ is the eigenvector
of purely imaginary $\alpha$-quaternion $I$. Here $X\in\ker(\id+I)$ or
$X\in\ker(\id-I)$, i.e. $X\in\ker\,q,\ |q|=0$. Inversely, let $X\in\ker\,q,\
|q|=0$. We have $q=a+b\tilde q;\quad a,b\in{\bf R},\ \ \tilde q$ is purely
imaginary $\alpha$-quaternion. Then $q(X)=0$, i.e. $\tilde q(X)=-\frac abX;
\frac ab=\pm 1$, thus $X$ is the eigenvector of purely imaginary $\alpha$
-quaternion $\tilde q$. It is easy to check that in this case $Orb(X)$ consists
of kernel elements of some isotropic $\alpha$-quaternions. Indeed, if $X$ is
eigenvector of purely imaginary $\alpha$-quaternion $I$, i.e. $IX=\pm X;\,
q=a\id+bI+cJ+dK\in{\bf H}_\alpha$ is any $\alpha$-quaternion, then
$\{\,\tilde q\in{\bf H}_\alpha\mid\tilde q(qX)=0\,\}=\{(u^2+v^2)x\,\id
+(\mp(u^2-v^2)x-2uvy)I+(-2uvx\pm u^2-v^2)y)J+(u^2+v^2)yK\mid x,y\in{\bf R}\,\}$
where $u=a+b,v=c-d$.
Therefore, we have proved
\procl{Theorem 2} The dimension of the linear representation orbit of
$\alpha$-quaternion algebra is equal to either 2 or 4. It is 2 iff the orbit
belongs to unification of kernels of isotropic $\alpha$-quaternions.\qed
\eprocl
\par
Like in Hermitian geometry we call 4-dimensional orbits of linear represen\-tation
of $\alpha$-quaternion algebra {\it holomorphic subspaces}. 2-dimensional orbits
of the representation are called ({\it real\/}) {\it spinor planes}.
\section{Generalized Almost Quaternionic Structures}
In this chapter we construct explicitly a fibre bundle of 4-dimensional
oriented pseudo-Riemannian manifold on the basis of the well-known notions
of self-dual and anti-self-dual forms on the manifold, the type fibre of the
bundle being isomorphic to the algebra of generalized quaternions. Such fibre
bundle  over manifolds of greater dimension is issential in the notion of a
generalized almost quaternionic structure introduced and investigated in the
chapter. We study in detail the so-called generalized almost quaternionic
structures of spinor type, generalized almost quaternionic structures with
parallelizable structural bundle being their particular case. The interrelation
of the types of generalized almost quaternionic structures is also investigated.
\subsection{Features of 4-dimensional pseudo-Riemannian mani\-folds geometry}
Let $(M^n,g)$ be an $n$-dimensional oriented pseudo-Riemannian manifold,
${\sf X}(M)$ be the module of smooth vector fields on $M$ over ring $C^\infty(M)$
of smooth functions on $M,\ \ TM=\cup_{p\,\in M}T_p(M)$ be a tangent
bundle over $M,\ \ {\cal T}(M)=\oplus_{r,s=0}^\infty{\cal T}^s_r(M)$
be a tensor algebra, and $\Lambda(M)=\oplus_{r=0}^\infty\Lambda^r(M)$
be a Grassmann algebra of differential forms of manifold $M,\ \{\,E\,\}$ be
a module of smooth sections of fibre bundle $(E,M,\pi)$ over $M$. All manifolds,
tensor fields and other objects are considered smooth and of class $C^{\,\infty}$.
\par
Recall [1] that Hodge operator $*:\Lambda(M)\rightarrow\Lambda(M)$ is defined
by the sequence of (necessarily unique) operators $*:\Lambda^{r}(M)\rightarrow
\Lambda^{n-r}(M);\quad r=0,\ldots,n$, such that $\omega\land(*\vartheta)=\la\omega,
\vartheta\ra\eta_{g}$, where $\omega,\vartheta\in\Lambda^r(M),\ \eta_g$ is
volume form on $M,\ \la\cdot,\cdot\ra$ is metric on $\Lambda^r(M)$ induced by metric
$g=\la\cdot,\cdot\ra$. We know that $*^2=(-1)^{r(n-r)+s}\id$, where $s$ is the negative
inertia index (or, the index) of metric $g$.
\par
Let, in particular, $n=4,r=2.$ Then operator $*:\Lambda^2(M)\rightarrow
\Lambda^2(M)$ is involutory, i.e. $*^2=\id$, iff $s$ is even, i.e. $s=0$ or 2
(we take $n-s\ge s$). In this case $\Lambda^2(M)=\Lambda^+(M)\oplus\Lambda^-(M)$,
where $\Lambda^+(M)$ and $\Lambda^-(M)$ are eigensubmodules of endomorphism "$*$",
corresponding to eigenvalues 1 and (-1), respectively. Their elements are
called,respectively, {\it self-dual} and {\it anti-self-dual} forms [1].
\par
It is known that giving tensor $t$ of the type $(r,s)$ on smooth manifold is
equivalent to giving a system of smooth functions $\{\,t^{j_1...j_s}{}_{i_1...i_r}\}$
on the space of the frame bundle over the manifold, the functions satisfying
the known system of differential equations [28]. The functions are called
{\it components} of tensor $t$. Denote by $\{\,\eta_{ijkl}\,\}$ components of tensor
$\eta_g$ in a positively oriented frame. Then, if $\{\,\omega _{ij}\,\}$ are
components of 2-form $\omega\in\Lambda^2(M)$, then $(*\omega)_{ij}=\frac12\eta_{ijkl}\,\omega^{kl}$,
i.e. $*\omega_{ij}=\frac12 g^{km}\,g^{lr}\,\eta_{ijkl}\,\omega_{mr}$,
where $\{\,g^{ij}\,\}$ are contravariant components of metric tensor. Thus,
\begin{eqnarray}
\omega\in\Lambda^+(M)\quad &\iff\quad\omega_{ij}
=\frac12\eta_{ijkl}\,g^{km}\,g^{lr}\,\omega_{mr};\nonumber\\
\omega\in\Lambda^-(M)\quad &\iff\quad\omega_{ij}=-
\frac12\eta_{ijkl}\,g^{km}\,g^{lr}\,\omega_{mr}.\end{eqnarray}
Giving a pseudo-Riemannian structure $g$ on an oriented manifold is equivalent
to giving a $G$-structure having the structural group $G=SO(n,s;{\bf R})$.
The elements of the $G$-structure are called {\it oriented orthonormal frames}.
In case $n=4,s=0$ or 2 we suppose that vectors of the frames are numerated
so that on the space of the $G$-structure $(g_{ij})=diag(1,-\a,-\a,1)$,
where $\a=s-1$. Then correlations (7) have the following form on the space
of the $G$-structure:
\begin{equation}
\omega\in\Lambda^\pm (M)\quad\iff\quad\omega_{\hat i\hat j}=
\pm\e(i,j)\,\omega_{ij},
\end{equation}
where $\e(i,j)=g_{ii}\,g_{jj}\,,\ (\hat i,\hat j)$ is a pair complementing
the pair $(i,j)$ up to even permutation of indices (1,2,3,4).
\procl{Theorem 3} Giving self-dual form $\omega$ on 4-dimensional oriented
pseudo-Rieman\-nian manifold $(M,g)$ of index $s=0$ or $2$ is equivalent to
giving an endomorphism $J$ of module ${\sf X}(M)$ connected with the form
$\omega$ by the identity
\begin{equation}
\omega(X,Y)=\la X,JY\ra;\qquad X,Y\in{\sf X}(M),
\end{equation}
and having the properties:
\begin{eqnarray}
1.\ \la JX,Y\ra+\la X,JY\ra=0;\qquad X,Y\in{\sf X}(M);\nonumber\\
2.\ J^2=-\lambda^2\id;\qquad\lambda^2=-\frac14\tr(J^2)=\frac12\|\omega\|^2;
\end{eqnarray}
\par{\rm3}.\ Orientation of manifold $M$ defined by the form $\omega\land\omega$
coincides with the orientation fixed on $M$ iff $\ \lambda\in{\bf R}\setminus
\{0\}$.
\eprocl
\demo{Proof} Let $\omega\in\Lambda^+(M)$. In view of (8) it implies that
on the space of $G$-structure
\begin{equation}
(\omega_{ij})=\left(\begin{array}{rrrr}
0&x&y&z\\-x&0&z&\alpha y\\-y&-z&0&-\alpha x\\
-z&-\alpha y&\alpha x&0\end{array}\right)
\end{equation}
Raising index of the form $\omega$ we get the matrix of endomorphism $J$
components on the space of $G$-structure:
\begin{equation}
(J^i{}_j)=\left(\begin{array}{rrrr}
0&x&y&z\\ \alpha x&0&-\alpha z&-y\\ \alpha y&\alpha z&0&x\\
-z&-\alpha y&\alpha x&0\end{array}\right)
\end{equation}
Note that metrix $g$ of index $s$ induces on module $\Lambda^+(M)$ a metric
of the same index. Indeed,
$$
\|\omega\|^2=\frac12\omega_{jk}\,\omega^{jk}=\frac12\sum_{i,j=1}^4\epsilon(j,k)(\omega_{jk})^2,
$$
where $\epsilon(j,k)=g_{jj}g_{kk}$, and in view of (11),\quad
$\|\omega\|^{2}=2(z^2-\alpha x^2-\alpha y^2)$. Now, it is clear that $(10_{1})$
simply means skew-symmetry of the form $\omega,\  (10_2)$ is checked directly
by (12), and $(10_3)$ follows from the definition of Hodge operator.
\par
Inversely, let $J$ be an endomorphism of module ${\sf X}(M)$ having properties
(10) and $\omega$ be the 2-form defined by (9). Then on the space of $G$-structure
$$
(\omega_{ij})=\left(\begin{array}{rrrr}
0&x&y&z\\-x&0&u&v\\-y&-u&0&w\\-z&-v&-w&0\end{array}\right);\quad
(J^i_j)=\left(\begin{array}{rrrr}
0&x&y&z\\\alpha x&0&-\alpha y&\alpha v\\
\alpha y&\alpha x&0&-\alpha w\\-z&-v&-w&0\end{array}\right)
$$
Hence $J^2=-\lambda^2\id$ iff
$$\begin{array}{llll}
&1)\ uy=\alpha vz;\qquad &2)\ ux=-\alpha wz;\qquad &3)\ ux=-wy;\\
&4)\ xy=-uw;\qquad &5)\ xz=-\alpha vw;\qquad &6)\ yz=\alpha uv;\end{array}$$
$$
\begin{array}{l}7)\ \alpha(x^2+y^2)-z^2=\alpha(x^2+v^2)-u^2=\\\quad\ \alpha
(y^2+w^2)-u^2=\alpha(v^2+w^2)-z^2=-\lambda^2.\end{array}$$
It is easy to see that the system is equivalent to the equations:
\begin{equation}
z=\pm u;\quad y=\pm\alpha v;\quad x=\mp\alpha w;\quad\lambda^2=z^2-\alpha x^2-\alpha y^2.
\end{equation}
On the other hand, $(\omega\land\omega)_{1234}=2(\omega_{12}\,\omega_{34}
+\omega_{13}\,\omega_{42}+\omega_{14}\,\omega_{23})=2(xw-yv+zu)$, and in view
of (13), $(\omega\land\omega)_{1234} =\mp 2(\alpha x^2+\alpha y^2-z^2)=2\lambda^2$.
Thus, (10) hold iff $z=u,\ y=\alpha v,\ x=-\alpha w$, i.e.
$\omega \in\Lambda^+(M)$.\qed
\edemo
\procl{Proposition 7} Let $\omega,\vartheta\in\Lambda^+(M); J_1,J_2$ be
corresponding endomorphisms of module ${\sf X}(M)$. Then $\omega\bot\vartheta$
in the induced metric of Grassmann algebra iff $J_1\circ J_2+J_2\circ J_1=0$.
\eprocl
\demo{Proof} In view of $(10_2),J^2=-\frac12\|\omega\|^2\id\quad(\omega\in\Lambda^+(M))$.
Polarize the equation. Let $\omega,\vartheta\in\Lambda^+(M),\ I$ and $J$ be
their corresponding endomorphisms. Then $(I+J)^2=-\frac12\|\omega+\vartheta\|^2\id$.
Opening the brackets and by (10) we have
$$I\circ J+J\circ I=-\la\omega,\vartheta\ra\id,$$
and the assertion being proved immediately follows.\qed
\edemo
Identify self-dual forms on $M$ with corresponding endomorphisms of module
${\sf X}(M)$ and denote by ${\cal S}c\,(M)\subset{\cal T}^1_1(M)$ a subbundle
of the fibre bundle of scalar endomorhisms on manifold $M$.
\procl{Theorem 4} Let $M$ be a 4-dimensional oriented pseudo-Riemannian
manifold of index $0$ or $2$. Then with the given identification fibre bundle
${\sf O}={\cal S}c\,(M)\oplus\Lambda^+(M)$ over $M$ is a subbundle of fibre bundle
${\cal T}^1_1(M)$,its standard fibre being isomorphic to algebra ${\bf H}_\alpha$.
\eprocl
\demo{Proof} It follows from the above that in every point $p\in M$ the fibre
of bundle $\sf O$ is a subalgebra of algebra $\End T_p(M)$. Choose in space
$\Lambda^+_p(M)$ an orthogonal basis $\{\,\omega_1,\omega_2,\omega_3\,\}$ whose
elements have the norms $\sqrt{-2\alpha},\sqrt{-2\alpha},\sqrt{2\alpha},$
respec\-tively. Let $J_1$ and $J_2$ be endomorphisms corresponding to the forms
$\omega_1$ and $\omega_2$, respectively. In view of Proposition 7
$\ J_1\circ J_2+J_2\circ J_1=0.$ Besides, by $(10_2),\ (J_1)^2=(J_2)^2=\alpha\id$
and thus, the given subalgebra is isomorphic to algebra ${\bf H}_\alpha$.\qed
\edemo
\remark{Remark} In the similar way, it is checked that the anti-self-dual form
$\omega\in\Lambda^-(M)$  on the space of $G$-structure is defined by the matrix
of components
$$(\omega_{ij})=\left(\begin{array}{rrrr}
0&x&y&z\\-x&0&-z&-\alpha y\\-y&z&0&\alpha x\\-z&\alpha y&-\alpha x&0
\end{array}\right)$$
Here, the properties of self-dual forms on $M$ formulated in Theorems 3 and 4,
as well as in Proposition 7, hold for anti-self-dual forms, too, except $(10_3)$
that is here substituted by:
\par
{\it "Orientation of manifold M, defined by the form} $\omega\land\omega,$
{\it where} $\omega\in\Lambda^-(M)$ {\it is inverse to the orientation fixed on M
iff $\lambda\in{\bf R}\setminus\{0\}$.{\it"}
\eremark
\definition {Definition 3} Fibre bundle $\sf O$ over 4-dimensional oriented
pseudo-Rieman\-nian manifold $M$ of index 0 or 2 constructed on the basis of
fibre bundle of self-dual (anti-self-dual, respectively) forms on $M$ in the
above  way is called {\it a canonical almost} $\alpha${\it-quaternionic
structure of the first}} ({\it second}, respectively) {\it kind} on the manifold.
\edefinition
\subsection{Generalized almost quaternionic structures}
\definition {Definition 4} We call {\it an almost} $\alpha${\it -quaternionic} $(AQ_\alpha-)$
{\it structure} on manifold $M$ a subbundle ${\sf Q}$ of fibre bundle
${\cal T}^1_1(M)$, whose standard fibre is algebra ${{\bf H}_\a}$. The fibre
bundle $\sf Q$ will be called {\it structural bundle} of $AQ_\alpha$-structure.
The manifold with the fixed $AQ_\alpha$ -structure is called {\it an almost}
$\alpha ${\it -quaternionic} $(AQ_\alpha-)$ {\it manifold}.
\edefinition
\par
When $\alpha=-1,$ the notion of $AQ_\alpha$-structure coincides with the
classical notion of almost quaternionic (in other terms, almost quaternal)
structure [1], [13]. When $\alpha=1, AQ_\alpha$-structure will be called
{\it an almost antiquaternionic}, or {\it an almost quaternionic structure
of hyperbolic type} (the term is also accepted).
\par
Let $\sf Q$ be an almost $\alpha$-quaternionic structure on $M$. Evidently,
$\sf Q$ is a 4-dimensional vector fibre bundle metrized by the canonical
metric of $\alpha$-quaternion algebra and associated with the principle fibre
bundle $B{\sf Q}=(B({\sf Q}),M,\pi,G_\alpha)$, its total space elements being
fours $\{\,\id,J_1,J_2,J_3\,\}$ of endomorphisms of tangent space $T_p(M)$
in an arbitrary point $p\in M$ satisfying the correlations $(J_1)^2=(J_2)^2=
\alpha\id,J_1\circ J_2+J_2\circ J_1=0,J_3=J_1\circ J_2$ and defining oriented
orthonormal frames in fibres of fibre bundle $\sf Q$, and its structure
group being Lie group $G_\alpha=SO(2-\alpha,1+\alpha;{\bf R})$. Sections of
fibre bundle $\sf Q$ will be called $\alpha${\it -quaternions on M.}
\example{Example 1} Any 4-dimensional oriented pseudo-Riemannian manifold
of index $s=0$ or 2 has, as we have seen, two canonical almost $\alpha$
-quaternionic srtuctures (almost quaternionic, when $s=0$, and almost
antiquaternionic, when $s=2$).
\eexample
\example {Example 2} A Riemannian $4n$-dimensional manifold is called
({\it locally}){\it hyper-Kaehler} if its (restricted) holonomy group is
contained in $Sp(n)$. The notion of hyper-Kaerher manifold was introduced
by  E.Calabi [30] who constructed the first non-trivial examples of such
manifolds. It is known [1]  that a hyper-Kaehler manifold can be characterized
as a Riemannian manifold admitting two anti-commuting complex structures
$\{\,I,J\,\}$  parallel in Riemannian connection. Evidently, the algebraic
shell of these endomorphisms defines the almost quater\-nionic structure on the
manifold, the structural bundle being parallelizable by endomorphisms
$\{\,\id,I,J,K\,\}$, where $K=I\circ J$.
\par
Geometry of hyper-Kaehler manifolds was studied by many authors, for example,
M.Berger [31], A.Beauville [19] and others. In particular, M.Berger proved
that any hyper-Kaehler manifold is Ricci-flat [31]. A homogeneous hyper-Kaehler
manifold is flat, so it is decomposed into the product of tori by Euclidean
spaces. It is well-known that a 4-dimensional Riemannian manifold is hyper-Kaehler
iff it is a Kaehler and Ricci-flat manifold [1]. At present examples of
compact hyper-Kaehler manifolds of any $4n$-dimension are known [19,20].
On the other hand, E.Calabi [30] showed that the space of cotangent bundle
of complex projective space ${\bf C}P^n$ admits a full hyper-Kaehler metric.
Generalizing  E.Calabi's constuction, N.Hitchin, A.Karlshede, U.Lindstr\"om
and M.Ro\u cek developed a new factorization method [1] based on symplectic
reduction of Marsden and Weinstein [32] allowing to build new examples of
non-compact but complete hyper-Kaehler manifolds. With this method, for instance,
an example of hyper-Kaehler manifold used for construction of gravitational
multi-instanton of Hibbons and Hoking [33] was built.
\eexample
\example{Example 3} $4n$-dimensional Riemannian manifold is called
{\it quaternionic-Kaehler} if its holonomy group is contained in
$Sp(n)\cdot Sp(1)$. It is known [1] that a quater\-nionic-Kaehler
manifold can be characterized as a Riemannian manifold having an almost
quaternionic structure, the structural bundle being invariant with respect to
parallel translations in Riemannian connection. Since $Sp(1)\cdot Sp(1)=
SO(4,{\bf R})$, any 4-dimensional oriented Riemannian manifold is
quaternionic-Kaehler. Geometry of quaternionic-Kaehler manifolds was also
studied by many authors, for example M.Berger [17], J.Wolf [22],
D.V.Alekse\-evskii [23] and other. In particular, M.Berger proved [17] that
a quaternionic-Kaehler manifold $M^{4n}\ (n>1)$ is Einsteinian, and it is
locally hyper-Kaehler iff it is Ricci-flat. J.Wolf [22] received a complete
classification of quaternionic-Kaehler symmetric Riemannian manifolds of
nonzero Ricci curvature: if such manifold has a positive Ricci curvature it
is compact and is one of the following spaces:
$$\begin{array}{lll}
&1)\ Sp(n+1)\,/\,Sp(n)\cdot Sp(1)={\bf H}P^n;\qquad
&2)\ SU(n+2)\,/\,S(U(n)\cdot U(2));\\
&3)\ SO(n+4)\,/\,SO(n)\cdot SO(4);\qquad
&4)\ G_2\,/\,SO(4);\\
&5)\ F_4\,/\,Sp(3)\cdot Sp(1);\qquad
&6)\ E_6\,/\,SU(6)\cdot Sp(1);\\
&7)\ E_7\,/\,Spin(3)\cdot Sp(1);\qquad
&8)\ E_8\,/\,E_7\cdot Sp(1).\end{array}$$
If it has a negative Ricci curvature it is non-compact and is one of the following spaces:
$$\begin{array}{lll}
&1)\ Sp(n,1)\,/\,Sp(n)\cdot Sp(1);\qquad
&2)\ SU(n,2)\,/\,S(U(n)\cdot U(2));\\
&3)\ SO(n,4)\,/\,SO(n)\cdot SO(4);\qquad
&4)\ G_2^2\,/\,SO(4);\\
&5)\ F_4^{-20}\,/\,Sp(3)\cdot Sp(1);\qquad
&6)\ E_6^2\,/\,SU(6)\cdot Sp(1);\\
&7)\ E_7^{-4}\,/\,Spin(3)\cdot Sp(1);\qquad
&8)\ E_8^{-24}\,/\,E_7\cdot Sp(1).\end{array}$$
Moreover, D.V.Alekseevskii [21] proved that any compact homogeneous
quater\-nionic-Kaehler manifold is a Riemannian symmetric space. The question
of the existence of non-symmetric quaternionic-Kaehler manifolds was posed
by S.Kobayashi and J.Eells [34] and positively solved by D.V.Alekseevskii
[23] who received a complete classification of quaternionic-Kaehler manifolds
admitting transitive solvable motion group.
\eexample
\example{Example 4} Almost antiquaternionic structures were studied,up to now,
mainly in the case of parallelizable structural bundle. V.F.Kirichenko [35]
showed that such structure is canonically defined, for example, on the space
of tangent bundle over the manifold of affine connection, as well as on the
manifold carrying three-web structure. In the works by M.A.Akivis and A.M.Shelekhov,
and many of their pupils, the three-web theory has been developed. In particular,
they received a large number of concrete examples of three-webs having a number
of remarkable properties, for example, Grassmanian three-webs [26]. In the next
chapter we shall see the main notions of tree-webs theory finding their
graceful expression in terms of generalized almost quaternionic structure
theory. The latter can be regarded, in particular, as a surprising and fruitful
three-web theory generalization.
\eexample
\subsection{Almost $\alpha$-quaternionic connections}
\definition{Definition 5} An affine connection on $AQ_\alpha$-manifold
is called {\it an almost  quater\-nionic}, or $AQ_{\alpha }${\it -connection},
if module $\sf Q$ of structural bundle sections is invariant with respect
to all parallel translations generated by the connection.
\edefinition
For example, by definition of quaternionic-Kaehler manifold, Riemannian
con\-nection on such manifold is a $AQ_\alpha$-connection.
\procl {Proposition 8} On $AQ_\alpha$-manifold $M$ there always exists
an $AQ_\alpha$-con\-nection.
\eprocl
\demo{Proof} Let $\{\,U_a\,\}_{a\in A}$ be the covering of manifold $M$ trivializing
bundle $\sf Q$, and let $\{\,\id,\un Ia,\un Ja,\un Ka\,\}$ be the section of
bundle $B{\bf Q}$ over $U_a$. We will see below (Theorem 9) that on manifold
$U_a$ there exists connection $\un\nabla a$ whose tensors $\un Ia,\un Ja,\un Ka$
are covariantly constant, in particular, $\un\nabla a$ is an $AQ_\alpha$-connection
on $U_a$. Let $\{\,\varphi_a\,\}_{a\in A}$ be a partitioning of unit dominated
by covering $\{\,U_a\,\}_{a\in A}$. Then, evidently, $\nabla=\sum_{a\in A}\varphi_a\un\nabla a$
is an $AQ_\alpha$-connection on $M$.\qed
\edemo
\procl{Proposition 9} Any $AQ_\alpha$-connection $\nabla$ on $AQ_\alpha$-manifold $M$
induces a metric connection in fibre bundle $\sf Q$.
\eprocl
\demo{Proof} Let $f\in C^\infty(M),\quad q\in\{{\sf Q}\},\quad X,Y\in{\sf X}(M)$.
Then
$$\begin{array}{rcl}
\nabla_X(fq)Y\eq\nabla(fq(Y))-fq(\nabla_XY)\\
{}\eq X(f)q(Y)+f\nabla_X(q)+fq(\nabla_XY)-fq(\nabla_XY)\\
{}\eq X(f)q(Y)+f\nabla(q)Y,\end{array}$$
and in view of arbitrary $Y\in{\sf X}(M),\ \nabla_X(fq)=X(f)q+f\nabla_X(q)$,
\ i.c. $\nabla$ is a linear connection in bundle $\sf Q$. Further, let
$q_1,q_2\in\{\,{\sf Q}\,\}$. Since a $AQ_\alpha$-connection generates
differentiation of algebra ${\cal T}(M)$ and the differentiation is
permutational with contractions, then
$$\forall X\in{\sf X}(M)\implies\nabla_X(q_1\,q_2)=\nabla_X(q_1)q_2+q_1\nabla_X(q_2).$$
Besides, $\nabla_X\,(\o q)=\o{\nabla_X\,(q)}\quad(q\in\{\,{\sf Q}\,\}).$ Indeed,
in our trivialization
$$q=b\id+c\un Ia+d\un Ja+f\un Ka;\qquad b,c,d,f\in C^\infty(U_a).$$
Since tensors $\un Ia,\un Ja,\un Ka,\id$ are parallel in connection
$\un\nabla a$,
$$\un\nabla a_X(q)=X(b)\id+X(c)\un Ia+X(d)\un Ja+X(f)\un Ka,$$
thus,
$$\nabla_X\,(q)=\sum_{a\in A}\varphi_a\{\,X(b)\id+X(c)\un Ia+X(d)\un Ja+X(f)\un Ka\,\}.$$
On the other hand, $\o q=b\id-c\un Ia-d\un Ja-f\un Ka$; and hence,
$$\nabla_X(\o q)=\sum_{a\in A}\varphi_a\{\,X(b)\id-X(c)\un Ia-X(d)\un Ja
-X(f)\un Ka\,\}=\o{\nabla_X(q)}.$$
It follows that
$$\nabla_X(\o q_1\,q_2)=\o{\nabla_X(q_1)}q_2+\o q_1\,\nabla_X(q_2);\quad
\nabla_X(\o q_2\,q_1)=\o{\nabla_X(q_2)}q_1+\o q_2\,\nabla_X(q_1).$$
Summarizing the equations elementwise and by definition of metric
$\la\cdot,\cdot\ra$ we get
$$\nabla_X(\la q_1,q_2\ra)=\la\nabla_X\,(q_1),q_2\ra+\la q_1,\nabla_X(q_2)\ra,$$
i.e. $\nabla$ is a metric connection in $\sf Q$.\qed
\edemo
\definition{Definition 6} An $AQ_\alpha$-manifold with fixed $AQ_\alpha$-connection
will be called {\it calib\-rated} and the fixed $AQ_\alpha$-connection will be
called {\it its calibration.}
\edefinition
Further on all considered $AQ_\alpha$-manifolds will be implied calibrated.
\par
Let $M$ be $AQ_\alpha$-manifold. In fibre bundle $\sf Q$ we distinguish
two natural subbundles -- subbundle $\sf T$ of purely imaginary
$\a$-quaternions $q$ satisfying the condition $q^2=\a$, or $\a(\im q)^2=1$,
and subbundle $\sf C$ of $\alpha$-quaternions $q$ satisfying the condition
$\alpha(\im q)^2>0$. We call them, respectively, {\it twistor} and {\it conformal
bundles} over $M$ associated to structure $\sf Q$. We call total fibre space
of $\sf T$ {\it a twistor space} over $M$. The terminology can be explained
by the fact that in case $M$ is a 4-dimensional oriented Riemannian manifold,
the twistor bundle over $M$ with respect to canonical almost quaternionic
structure coincides with the classical twistor Penrose bundle [1]. In case
$M$ is an arbitrary quaternionic manifold, that is an almost quaternionic
manifold, admitting almost quaternionic torsion-free connection [5], the
twistor bundle over $M$  coincides with the twistor bundle introduce by
S.Salamon [5] and L.Be\'rard-Bergery [6].
\par
Sections of fibre bundle $\sf T$ will be called {\it twistors} on $M$.
Evidently, any twistor is either an almost complex structure or a structure of
almost product on $M$ when $\alpha=-1$ or 1, respectively. Note that in view
of Proposition 9 an $AQ_\alpha$-connection on $M$ induces a connection in
conformal bundle $\sf C$.
\subsection{Generalized almost quaternionic structure of spinor type}
\definition {Definition 7} An almost $\alpha$-quaternionic structure on manifold
$M$ is called {\it a structure of spinor type} if $M$ admits a twistor.\edefinition
\definition {Definition 8} An almost $\alpha$-quaternionic structure on $M$
is called {\it an} $\pi AQ_\alpha${\it -struc\-ture} if $M$ admits a pair of
anticommuting twistors.\edefinition
\par
The above condition is, evidently, equivalent to parallelizability of structural
bundle: if $\{\,I,J\,\}$ is a pair of anticommuting twistors on $M$, then
endomorphisms $\{\,\id,I,J,K\,\}$,  where $K=I\circ J$, define the parallelizma
of structural bundle.
\par
For example, any hyper-Kaehler manifold, as well as the manifold carrying
a three-web structure  is a $\pi AQ_\alpha$-manifold and, certainly, a
$AQ_\alpha$-manifold of spinor type. On the other hand, by [1], quaternionic
projective space ${\bf H}P^n$, being a quaternionic-Kaehler manifold, does not
admit an almost complex structure and, thus, is not an $AQ_\alpha$-manifold of
spinor type. By similar reason, 4-dimentional sphere $S^4$, being a
quaternionic-Kaehler manifold, is not an $AQ_\alpha$-manifold of spinor type.
\procl{Theorem 5} The space of conformal bundle over $AQ_\alpha$-manifold
$M$ admits a natural structure of an $AQ_\alpha$-manifold of spinor type.
\eprocl
\demo{Proof} Let $\vartheta$ be the connection form in fibre bundle $\sf C$
induced by calibration $\nabla,q\in{\sf C},q_1\in{\sf C}_{\pi(q)}(M),X\in T_q({\sf G})$.
Assume
\begin{equation}
Q_1(X)=\tau^{-1}\circ L_{q_1}\circ\tau\circ\vartheta(X)+\pi ^{-1}_*\circ q_1\circ\pi_*(X),
\end{equation}
where $\pi^{-1}_*$ is the operator of horizontal lifting from $T_{\pi(q)}(M)$
into horizontal distribution of connection $\vartheta,\tau$ is the operator
of identification of vertical subspace at point $q$ with fibre
${\sf C}_{\pi(q)}(M),L_{q_1}$ is the left shift in the fibre on $\alpha$-quaternion
$q_1$. Evidently, the family of operators constructed in this way defines an
$AQ_\alpha$-structure on the space of bundle $\sf C$, and endomorphism $J$
generated by operators $J_q=|\im q|^{-1}\cdot\im q$ at every point $q\in{\sf C}$
defines the twistor on space $\sf C$.\qed
\edemo
Denote the constructed $AQ_\a$-structure on manifold $\sf C$ by ${\sf Q}_c$
and call it {\it conformal lifting} of the initial $AQ_\alpha$-structure $\sf Q$
Theorem 5 reduces geometry of an arbitrary $AQ_\alpha$-structure $\sf Q$
to the geometry of a certain $AQ_\alpha$-structure of spinor type --- conformal
lifting of structure $\sf Q$.
\definition {Definition 9} An almost $\alpha$-quaternionic structure $\sf Q$
on manifold $M$ is called {\it integrable} if there exists atlas
$\{\,U_a,\varphi_a\,\}_{a_\in A}$ of manifold $M$ trivializing fibre bundle
$B(\sf Q)$, the sections
$\{\,I_a,J_a\,\}_{a\in A}$ defining trivialization over $U_a,\quad a\in A$
being given in the natural basis by constant matrices.
\edefinition
>From (4) immediately follows
\procl{Proposition 10} Integrability of $AQ_\alpha$-structure is equivalent
to integrability of its conformal lifting.\qed
\eprocl
Let $\sf Q$ be an almost $\alpha$-quaternionic structure of spinor type on
manifold $M,\,J$ be a fixed twistor on $M$. In fibre bundle $\sf Q$ there is
naturally introduced the structure of 2-dimensional Hermitian vector fibre
bundle $\cal S$ whose standard fibre is $\alpha$-spinvector space ${\cal S}_0$
and Hermitian metric is a Hermitian form $\la\la X,Y\ra\ra=\la X,Y\ra+j\la X,J^3Y\ra$.
We call the fibre bundle  {\it the bundle of} $\alpha${\it-spinvectors over M.}
It is associated to the principal bundle $B{\cal S}=(B({\cal S}),M,\pi,SCo(2,
{\bf K}_\alpha))$ whose structural group is Lie group ${\bf H}^*_\alpha=SCo(2,
{\bf K}_\alpha)=SU(2,{\bf K}_\alpha)\times{\bf R}^+=Sp_\alpha(1)\times
{\bf R}^+$ of non-isotropic really normalized $\alpha$-quaternions. The group
is regarded as one of special conformal transformations of $\alpha$-spinvector
space ${\cal S}_0$. The elements of total fiber space $B({\cal S})$ are
orthogonal frames of the fibres of fiber bundle ${\cal S},\ {\bf R}$-homothetic
to spinframes.
\procl{Proposition 11} There exists a one-to-one corespondence between
mani\-folds $B({\sf Q})$ and $B({\cal S})/{\bf R}^*$, where
${\bf R}^*$ is a multiplicative group of real numbers acting on the fibres of
fibre bundle as a group of homotheties.
\eprocl
\demo{Proof} Let $p=(m,\id,J_1,J_2,J_3)\in B(\sf Q)$ be a positively oriented
ortho\-normalized frame of space ${\sf Q}_m$. By Theorem 1 and Proposition 5 it
has a corresponding class $rU(1,{\bf K}_\alpha)$ of orthonormalized frames of
space ${\cal S}_m$, where $(J_k)=\sigma_k,\quad k=1,2,3$. Let $\tilde r\in
B({\cal S})\cap rU(1,{\bf K}_\a)$. Then, $\tilde r=rg$, where $g=diag(e^{ia},
e^{ia}),\ \det g=1$. Then, $e^{2ia}=1,$ that is $a=\pi k, g=\pm I$. Keeping
in mind that each orthonormalized frame $r$ has a corresponding class
$r{\bf R}^+$ of orthogonal frames with the same determining property, belonging
to the orbit of frame $r$ with respect to the structural group of fibre bundle
$B{\cal S}$, we get
$\{\,r\in B({\cal S})\mid(J_k)=\sigma_k,k=1,2,3\,\}=r{\bf R}^*$. The juxtaposition
$p\rightarrow r{\bf R}^*$ gives the desired correspondence.\qed
\edemo
\procl{Proposition 12} There exists an innerly defined $\Ad$-homomorphism
$S:B{\cal S}\rightarrow B{\sf Q}$ of principal fibre bundles juxtaposing to
the orthogonal frame $p=(m,\e_1,\e_2)\in B({\cal S})$ a positively
oriented orthonormalized frame $(m,\id,J_1,J_2,$ $J_3)\in B({\sf Q})$, such that
$(J_k)=\sigma_k;\quad k = 1,2,3.$
\eprocl
\demo{Proof} Let $p=(m,\e_1,\e_2)\in B({\cal S})$. Then
$S(pg)=(m,\id,J_1,J_2,J_3)\mid (J_k)_{pg}=\sigma_k;\quad k=1,2,3$. Note that
$(J_k)_{pg}=g^{-1}(J_k)_p\,g$, i.e. $(J_k)_p=g\sigma_k\,g^{-1}$. On the other
hand, $S(p)=(m,\id,\tilde J_1,\tilde J_2,\tilde J_3)\mid(\tilde J_k)_p=\sigma_k;
\quad k =1,2,3;\quad S(p)\Ad(g)=(m,\id,\hat J_1,\hat J_2,\hat J_3)\mid(\hat J_k)_p
=(\Ad(g)(\tilde J_k))_p=g(\tilde J_k)_pg^{-1}$. By comparison, we get
$[\hat J_k]_p=[J_k]_p$, i.e. $S(pg)=S(p)\,\Ad(g)$.\qed
\edemo
\procl{Lemma 1} Let $s:G_1\rightarrow G_2$ be a Lie group homomorphism.
Then $\forall g\in G_1\implies\Ad(sg)\circ s_*=s_*\circ\Ad(g)$, where
$s_*$ is differential of mapping $s$ in the Lie group unit.
\eprocl
\demo{Proof} Let $A_g$ be the inner Lie group automorphism generated by element
$g$. Then $(A_{s(g)}\circ s)h=s(g)s(h)s(g^{-1})=s(ghg^{-1})=(s\circ\Ad)h\quad
(h\in G_1)$. So, $A_{s(g)}\circ s=s\circ A_g$. We complete the proof passing
onto mapping differentials.\qed
\edemo
\procl{Proposition 13} Any connection in the principle fibre bundle
$(B({\sf Q}),M,\pi,$ $G_\a)$ induces a connection in the principle fibre bundle
$(B({\cal S}),M,\pi,SCo(2,{\bf K}_\alpha))$.
\eprocl
\demo{Proof} Let $\Theta$  be a connection form on $B({\sf Q})$. We define
the form $\vartheta=s_*(\Theta)$ on $B({\sf S})$. We have:
\begin{eqnarray*}
\vartheta_{pg}\eq s^{-1}_*\Theta_{S(p)g}\,s_*=s^{-1}_{*}\Ad(s(g^{-1}))\Theta_
{S(p)}\,s_*=\\\eq\Ad(g^{-1})s^{-1}_*\Theta_{S(p)}s_*=\Ad(g^{-1})\vartheta_p.
\end{eqnarray*}
Thus, $\vartheta$ is a connection form on $B({\cal S})$.\qed
\edemo
We call this connection {\it a spinor representation of initial connection}.
\par
Let $M$ be an almost $\a$-quaternionic manifold of spinor type. By Proposition 9,
its calibration induces a metric connection in fibre bundle $B{\sf Q}$.
Let $\vartheta$ be the form of spinor representation of the connection, $\cal H$
and $\cal V$ be horizontal and vertical distributions
on $B({\cal S})$, respectively. Fix $p\in B({\cal S}$). Let $S(p)=(m,\id,J_1,
J_2,J_3)$ be its corresponding positively oriented orthonormalized frame of
space ${\sf Q}_m$. Define endomorphisms ${\cal J}_k\in{\cal T}^1_1(B({\cal S}))$
by the formulas
\begin{equation}
({\cal J}_k)_p=\vartheta^{-1}\circ[J_k]\circ\vartheta+\pi^{-1}_*\circ J_k\circ
\pi_*;\qquad k = 1,2,3;
\end{equation}
where $\pi^{-1}_*$ is the operator of horizontal lifting from $T_m(M)$ to
${\cal H}_p,\,\vartheta^{-1}$ is the identification operator of Lie algebra
${\sf sp}_\alpha(1)\oplus t^1$ with ${\cal V}_p$. Evidently, the pair
$\{\,{\cal J}_1,{\cal J}_2\,\}$ defines an almost $\a$-quaternionic structure
with parallelizable structural bundle, that is a $\pi AQ_{\alpha}$-structure,
on manifold $B({\cal S})$.
\par
Denote ${\cal J}^V_k=\vartheta^{-1}\circ[J_k]\circ\vartheta$,
${\cal J}^H_k=\pi^{-1}_*\circ J_k\circ\pi_*$. Call the endomorphisms
{\it vertical} and {\it horisontal components of the structure}
$\{\,{\cal J}_1,{\cal J}_2\,\}$, respectively. Evidently, distributions $\cal V$
and $\cal H$ are invariant with respect to endomorphisms ${\cal J}^V_k$ and
${\cal J}^H_k$, respectively. Note that ${\cal J}^H_k=J^H_k$ is just horizontal
lifting of twistor $J_k\quad(k =1,2)$. In this notation (15) will have the
form ${\cal J}_k={\cal J}^V_k+{\cal J}^H_k$. Consider the dependence of the
endomorphisms on the right action of the structural group in fibre space. Let
$g\in{\bf H}^\times_\alpha$. Then
\begin{eqnarray*}
({\cal J}^V_k)_{pg}\eq\vartheta^{-1}_{pg}[J_k]_{pg}\vartheta_{pg}\\
\eq(\Ad(g^{-1})\vartheta_p)^{-1}[(J_k)_{pg}]\Ad(g^{-1})\vartheta_p\\
\eq\vartheta^{-1}_p(\Ad(g^{-1})^{-1})[(J_k)_{pg}]\Ad(g^{-1})\vartheta_p\\
\eq\vartheta^{-1}_pg\circ[(J_k)_{pg}]\circ g^{-1}\vartheta_p\\
\eq\vartheta^{-1}_p[(J_k)_p]\vartheta_p=({\cal J}^V_k)_p;\\
({\cal J}^H_k)_{pg}\eq(J^H_k)_{pg}=\pi^{-1}_*(s(g)J_k)\pi_*=(s(g)J_k)^H_p.
\end{eqnarray*}
Thus, we have proved
\procl{Proposition 14} With the given notation
$$1)\,({\cal J}^V_k)_{pg}=({\cal J}^V_k)_p;\quad 2)\,({\cal J}^H_k)_{pg}
=(s(g)J_k)^H_p.\qeds$$
\eprocl
Summarizing the above we see that on total space $B({\cal S})$ of the principle
fibre bundle $B{\cal S}$ over an almost $\a$-quaternionic manifold of spinor
type there is canonically induced a structure of an almost $\a$-quaternionic
manifold with parallelizable structural bundle generated by endomorphisms
$\{\,{\cal J}_1,{\cal J}_2\,\}$.
\definition{Definition 10} Manifold $B({\sf S})$, equipped with a
$\pi AQ_\a$-srtucture generated by endomorphisms $\{\,{\cal J}_1,{\cal J}_2\,
\}$ is called {\it a covering} $\cal S${\it-space} of $AQ_\a$-manifold $M$ of
spinor type, and a $\pi AQ_\a$-structure generated by the endomorphisms is
called {\it a covering} or {\it canonical} $\pi AQ_\alpha${\it-structure}.
\edefinition
By Proposition 14 we get that the natural projection $\pi$ of fibre space
$B{\cal S}$ generates "projecting" of endomorphisms $\{\,({\cal J}_k)_p\,\}$,
\,$p\in B({\cal S})$, of the covering $\pi AQ_\alpha$-structure into tensor
fibre space of the type (1,1) on $M$,  whose image coincides with the space of
structural bundle $\sf Q$ of the  initial $AQ_\alpha$-structure. Thus, the
covering $\pi AQ_\a$-structure is defined by the initial $AQ_\a$-structure
and, in its turm, defines it. We get
\procl{Theorem 6} Giving of almost $\a$-quaternionic structure of spinor
type on manifold is equivalent to giving of its covering $\pi AQ_\a$-structure
on the covering $\cal S$-space.\qed
\eprocl
\subsection{Spintensor algebra of generalized almost quaternionic manifold of spinor type}
Let $M$ be an almost $\alpha$-quaternionic manifold of spinor type,
$\cal S$ be a bundle of $\alpha$-spinvectors over $M$. By means of anti-automorphism
$\tau$ of complex conjugation in module ${\cal S}_0$\ (the standard fibre of
the bundle) we can introduce another structure of ${\bf K}_\a$-module, assuming
$z\circ X=\tau(z)\,X$. Denote  module ${\cal S}_0$ equipped with such structure by
$\o{\cal S}_0$ and call it {\it a module conjugated to module} ${\cal S}_0$.
We also introduce module ${\cal S}_0^*$, dual to module ${\cal S}_0$, and module
$\o{\cal S}{}_0^*$, dual to module $\o{\cal S}_0$. The section module of direct sum
${\cal S}\oplus\o{\cal S}\oplus{\cal S}^*\oplus\o{\cal S}{}^*$ of corresponding
bundles generates in a standard way a tensor algebra ${\cal S}(M)$ that we call
{\it a spinors algebra of manifold M}, its elements being called
$\alpha${\it-spinor on M.} The algebra has a natural 4-graduation:
$\alpha$-spinor of the type $\left[\begin{array}{cc}a&b\\c&d\end{array}\right]$
can be regarded as a linear mapping of the module
$$
(\otimes^a\,\{\,{\cal S}\,\})\otimes(\otimes^b\,\{\,\o{\cal S}\,\})\otimes
(\otimes^c\{\,{\cal S}^*\,\})\otimes(\otimes^d\,\{\,\o{\cal S}{}^*\,\})
$$
in $C^\infty\,(M)$. Finally, tensor product ${\cal S}(M)\otimes{\cal T}(M)$ is
called {\it a spintensor algebra of manifold} $M$ and denote by ${\cal ST}(M)$
its elements being called $\alpha${\it-spintensors on} $M$.
\par
Let $\vartheta$ be the form of spinor representation of metric connection in
fibre bundle $B{\sf Q}$ induced by calibration of manifold $M$. We construct
the mapping $f$ of vertical right-invariant vector field module on manifold
$B({\sf Q})$ onto module of $\alpha$-quaternions on $M$, assuming
$(f(X))_\pi(p)=S(p)\,\vartheta_p(X_p)$, where $S(p):{\bf H}_\a\rightarrow
{\sf Q}_{\pi(p)}$ is a canonical mapping defined by frame $S(p)$. We show the
correctness of the definition in the sense of independence on the choice of
frame $p\in B({\cal S})$.
Let $\pi(p)=\pi(p');\,p'=pg,\,g\in SCo(2,{\bf K}_\alpha)$. Then
\begin{eqnarray*}
S(p')\,\vartheta_{p'}(X_{p'})\eq S(pg)\,\vartheta_{pg}(X_{pg})\\\eq S(p)\,
\Ad(g)\Ad(g^{-1})\,\vartheta_p(X_p)\\\eq S(p)\,\vartheta_p(X_p)=(f(X))_{\pi(p)},
\end{eqnarray*}
i.e. the definition is correct.
\par
Thus, the mapping $f$ identifies vertical right-invariant vector fields on
manifold $B({\cal S})$ with $\alpha$-quaternions on manifold $M$. On the other
hand, if $X$ is an arbitrary right-invariant vector field on $B({\cal S})$,
then $X=X^V+X^H$ and in view of canonical identification of right-invariant
horizontal vector fields on $B({\cal S})$ with elements from ${\sf X}(M)$ we get:
\procl{Theorem 7} Module ${\sf X}^R(B({\cal S}))$ of right-invariant vector
fields on the covering $\cal S$-space of $AQ_\alpha$-manifold of spinor type
is canonically identified with the direct sum of modules $\{\,{\sf Q}\,\}
\oplus{\sf X}(M)$.\qed
\eprocl
Note that module $\{\,{\sf Q}\,\}$ is naturally identified with a submodule of
spinor algebra of manifold $M$ generated of Hermitian form $\la\la\cdot,\cdot\ra\ra$
and the submodule skew-Hermitian forms on $\{\,{\sf Q}\,\}$. Indeed,
let $q$ be an arbitrary $\alpha$-quaternion on $M,\,q=(\re q)\id+\im q$. Assume
\begin{eqnarray*}
Q(X,Y)\eq\la\la X,q(Y)\ra\ra\\\eq\la\la X,(\re q)Y\ra\ra+\la\la X,(\im q)Y\ra
\ra\\\eq\la\la X,Y\ra\ra\,\re q+\la\la X,(\im q)Y\ra\ra;\qquad X,Y\in\{\,
{\cal S}\,\}
\end{eqnarray*}
Evidently,
\begin{eqnarray*}
Q(X,Y)\eq\la\la X,(\im q)Y\ra\ra=\re(\o X\cdot\im q\cdot Y)=\re(\o{\o Y\cdot
\o{\im q}\cdot X})\\
\eq\o{\re(\o Y\cdot\o{\im q}\cdot X)}=-\o{\re \o Y\cdot\im q\cdot X)}
=-\o{\la\la Y,(\im q)X\ra\ra}\\\eq-\o{Q(Y,X)},
\end{eqnarray*}
i.e. $Q$ is a skew-Hermitian form. In view of this identification the direct
sum $\{\,{\sf Q}\,\}\otimes{\sf X}(M)$ generates a subalgebra of spintensor
algebra of manifold $M$. Thus, by Theorem 7, we get
\procl{Theorem 8} There exists a canonical non-graduated monomorphism of
tensor algebra of right-invariant tensor fields on the covering $\cal S$-space
of an $AQ_\alpha$-manifold of spinor type into the spintensor algebra of the
manifold.\qed
\eprocl
The monomorphism juxtaposes the right-invariant tensor $t\in{\cal T}^R(B({\cal S}))$
to a defined set of spintensors $\{\un tk\,\}$ on $M$, where
$\un tk(X_1,\ldots,X_m)=t(\psi X_1,\ldots$ $\ldots,\psi X_m)$, where $\psi X_j$ is
either $f^{-1}(X_j)\in{\cal V}$ or ${\cal V}^*$, either $\pi^{-1}_*(X_j)\in{\cal H}$
or ${\cal H}^*,\quad j = 1,\ldots,m$,  depending on whether there is 0 or 1 in the
$j$-th place of the binary representation of $k$. The set $\{\un tk\,\}$
of spintensors defining the right-invariant tensor $t$ is called {\it the spectre}
of tensor $t$.
\section{Generalized Almost Quaternionic Manifolds with Parallelizable
Structural Bundle}
The notion of an almost quaternionic structure was, evidently, first considered
by P.Libermann [8] and further studied by M.Obata [10], S.Ishihara [11] and
many other authors. The almost quaternionic structures studied in mentioned
papers were defined by giving of two anticommuting almost complex structures
and, thus, had a parallelizable structural bundle. The same feature was characteristic
of almost antiquaternionic structures generated by geometry of tangent bundle,
three-web structure, as well as by almost Hermitian structure of hyperbolic
type on Riemannian manifold [35]. The present chapter is devoted to the study
of the structures from a single view point as $\pq$-structures. This class of
$AQ_\a$-structures is of importance since (as we saw in the previous chapter)
the study of any $AQ_\a$-structure is generally reduced to the study of its
covering $\pq$-structure. The main result of the present chapter is having
proved that $\pq$-structure on a manifold innerly generates a unique
$AQ_\a$-connection of a special type that we call {\it canonical the connection},
in particular, generalizes the Chern connection widely used in three-web theory [36].
A formula for global definition of canonical connection is received. Criteria
of integrability of $\pq$-structure and its structural endomorphisms are found
in terms of torsion and curvature tensors of canonical connection. We also
introduce the notions of isoclinic distribution and isoclinic $\pq$-structure,
as well as isoclinic-geodesic $\pq$-structure, that generalizes the corresponding
notions of three-web theory [37]. Finally, we find a criterion of semi-holonomic
$\pq$-structure being isoclinic generalizing the well-known criterion of
M.A.Akivis in three-web theory [37].
\subsection{Canonical connection of $\pq$-structure}
Recall that an almost $\a$-quater\-nionic structure with parallelizable
structural bundle is uniquely defined by a pair of anticommuting twistors.
Consider them fixed and assume the following definition equivalent to Definition 8:
\definition {Definition 11} A pair $\{\,I,J\,\}$ of endomorphisms of module
${\sf X}(M)$ for which
$$
1)\,I^2=J^2=\a\id;\quad 2)\,I\circ J+J\circ I=0\qquad(\a=\pm1),
$$
is called $\pq${\it-structure on manifold M}.
\edefinition
The mentioned definitions are equivalent since the algebraic shell of endomor\-phisms
$I$ and $J$  defines subbundle $\sf Q$ of a fibre bundle of tensors of the
type (1,1) on $M$ with standard fibre ${\bf H}_\a$, that is an $AQ_\a$-structure
on $M,\ I$ and $J$ being a pair of anticommuting twistors with respect to the
structure. We call $I,J$ and $K=I\circ J$ {\it structural endomorphisms\/} of
$\pq$-structure.
\procl{Theorem 9} Any $\pq$-manifold $\{\,M,I,J\,\}$ admits a unique
connection $\nabla$ having the following properties:
\begin{equation}
1)\,\nabla I=0;\quad2)\,\nabla J=0;\quad3)\,S(IX,Y)-S(X,IY)=0;
\end{equation}
where $\ S(X,Y)=\n XY-\n YX-[X,Y]\ $ is torsion tensor of the connection $\nabla$,
$[\cdot,\cdot]$ are Lie brackets, $X,Y\in{\sf X}(M)$.
\eprocl
\demo{Proof} In ${\bf K}_\a\otimes C^\infty\,(M)$-module ${\bf K}_\a\otimes
{\sf X}(M)$ there are naturally defined four pair-wise conjugate submodules --- the
eigendistributions $D^i_I,D^{-i}_I, D^i_J,D^{-i}_J$ of structural endomorphisms
with eigenvalues $i,-i,i,-i$, respectively. Besides, there are natural\-ly
defined endomorphisms $V=\frac 12(\id+iI^3)$ and $H=\frac12(\id-iI^3)$ that
are projectors complementary to each other on distributions $D^i_I$ and
$D^{-i}_I$, respectively. Let $\ X,Y\in{\sf X}(M),\,X=VX+HX,\,Y=VY+HY$. Denote
$VX=a,\,VY=b,\,HX=u,\,HY=v$. Then $X=a+u,\,Y=b+v$. Assume that the desired
connection exists. Then
\begin{eqnarray}
\n X(IY)\eq\n X(I)Y+I\n XY=I\n XY;\\
\n X(JY)\eq\n X(J)Y+J\n XY=J\n XY.
\end{eqnarray}
Further, $iS(a,u)=S(ia,u)=S(Ia,u)=S(a,Iu)=S(a,-iu)=iS(a,u)$, hence $S(a,u)=0$,
and thus, $[a,u]=\n au-\n ua-S(a,u)=\n au-\n ua$. But by (17),
$\n au=\n a(HX)=H\n aX\in D^{-i}_I,\ \n ua=\n u(VY)=V\n uX\in D^i_I$. Then,
$$\n au=H[a,u],\quad \n ua=-V[a,u]=V[u,a].$$
But then in view of (18),
$$\begin{array}{rclcl}
\n ab\eq\n a(J \tilde u)=J\,\n a\tilde u=JH[a,\tilde u]=\a JH[a,Jb]\\
\eq\a VJ[a,Jb]=VJ^3\,[a,Jb];\\
\n uv\eq\n u(J\tilde b)=J\n u\tilde b=JV[u,\tilde b]=\a JV[u,Jv]\\
\eq\a HJ[U,jV]=HJ^3[u,Jv].
\end{array}$$
Therefore,
$$\n XY=\n ab+\n av+\n ub+\n uv=VJ^3[a,Jb]+H[a,v]+V[u,b]+HJ^3[u,Jv],$$
or in view of the notations,
\begin{equation}
\n XY=V\{\,[HX,VY]+J^3\,[VX,JVY]\,\}+H\{\,[VX,HY]+J^3\,[HX,JHY].
\end{equation}
In view of definitions $V$ and $H$,\ (19) we can write explicitly the following:
\begin{eqnarray}
\n XY\eq\frac 14\{\,[X,Y]-\a[IX,IY]+\a J[X,JY]-J[IX,IJY]-\nonumber\\
&{}&-\a I[IX,Y]+\a I[X,IY]-IJ[X,IJY]+IJ[IX,JY]\,\}.
\end{eqnarray}
Inversely, immediate checking show that $\nabla$ is affine connection on $M$
having the required properties.\qed
\edemo
We call the constructed connection {\it a canonical connection of a}
$\pq$-{\it structure}. Identities $16_1$ and $(16_2)$ show that when
$\a=-1$ connection $\nabla$ is almost quaterni\-onic in the sense of M.Obata
[10]. On the other hand, considering the almost antiquaternionic structure joined
to a three-web [35] easily shows that invariant connection joined to a three-web [36]
later called  {\it Chern connection\/} has the properties of the canonical
connection of adjointed $\pq$-structure and, thus, coincides with it.
\procl{Theorem 10} $\pq$-structure is integrable iff torsion and curvature
tensors of its canonical connection are zero.
\eprocl
\demo{Proof} Integrability of $\pq$-structure $\{\,I,J\,\}$ means the existence
of an atlas of manifold $M^n$ carrying the structure, tensors $I$ and $J$
having constant coordinates in every map of the atlas. Fix map $(U,\varphi)$
of such atlas with local coordinates $(x^1,\ldots,x^n)$ and let $\{\,e_k\,\}$
be the natural map basis, $e_k=\frac{\partial\ \ }{\partial x^k},\quad k=1,\ldots,n=
\dim M$. Here $\varphi:U\rightarrow V\subset R^n$ is the mapping of the map.
Any vector $X\in{\bf R}^n$ generates an one-parametric group of diffeomorphisms
$\{\,T_{t}\,\},\quad T_{t}(Y)=Y+tX$ of the space ${\bf R}^n$ and, respectively,
a local one-parameter group of diffeomorphisms of domain $V$ that we denote by
the same symbol. Evidently, $T_t(\tilde e_k)=(T_t)_*\,\tilde e_k=\tilde e_k$,
where $\tilde e_k=\varphi_*(e_k),\ f_*$ is differential of mapping $f$. Accordingly,
in domain $U$ there acts a local one-parametric group of diffeomorphisms
$\{\,F_t\,\}$ where $F_t=\varphi^{-1}\circ T_t\circ\varphi$. Evidently,
$(F_t)_*\,e_k=(\varphi^{-1}\circ T_t)_*\,\tilde e_k=(\varphi^{-1})_*\,\tilde e_k=e_k$.
Thus, ${\cal L}_{(\varphi^{-1})_*\,X}\,(e_k)=0,$ where $\cal L$ is the operator of
Lie differentiation. In particular, ${\cal L}_{e_i}(e_j)=0$, that is
\begin{equation}
[e_i,e_j]=0;\qquad i,j=1,\ldots,n.
\end{equation}
Without the loss of generality it can be considered that the given basis has
the form $\{\,e_1,\ldots,e_m,Ie_1,\ldots,Ie_m\,\}$. Introduce the complex
basis $\{\,\e_1,\ldots,\e_m,$ $\e_{\hat 1},\ldots,\e_{\hat m}\,\}$, where
$\e_a=\frac 12(\id+iI^3)e_a$ and $\e_{\hat a}=\frac 12(\id-iI^3)e_a,
\quad a=1,\ldots, m$. Evidently, $\{\,\e_a\,\}$ is the basis of distribution
$D^i_I|U$. From (21) we have that $[\e_a,\e_{\hat b}]=0.$
But then $\nabla_{\e_a}{\e_{\hat b}}=H[\e_a,\e_{\hat b}]=0$, hence in the above
notations $\n au=\nabla_{{a^b}\e_b}(u^{\hat c}\,\e_{\hat c})
=a^b\,\nabla_{\e_b}(u^{\hat c}\,\e_{\hat c})=a^b\,\e_b\,(u^{\hat c})\e_{\hat c}
=a(u^{\hat c})\e_{\hat c}$. Similarly, $\n ua=u(a^c)\e_c$. Further,
$\nabla_{\e_a}{\e_b}=VJ^3[\e_a,J\e_b]=0,\,\nabla_{\e_{\hat a}}{\e_{\hat b}}
=HJ^3[\e_{\hat a},J\e_{\hat b}]=0$, hence, $\n ab=a(b^c)\e_c;\,\n uv
=u(v^{\hat c})\e_{\hat c}$, and thus, $\n XY=X(Y^a)\e_a+X(Y^{\hat a})
\e_{\hat a}=X(Y^k)\e_k,\quad k=1,\dots,n$. Consequently,
\begin{eqnarray*}
S(X,Y)\eq\n XY-\n YX-[X,Y]=X(Y^i)\e_i-Y(X^i)\e_i-[X,Y]^i\e_i=0;\\
R(X,Y)Z\eq\nabla_X\n YZ-\nabla_Y\n XZ-\nabla_{[X,Y]}\,Z\\
\eq X(Y(Z^k))\e_k-Y(X(Z^k))\e_k-[X,Y](Z^k)\e_k=0.
\end{eqnarray*}
Inversely, if torsion and curvature tensors of canonical connection are zero,
then there exists an atlas where Christoffel connection  coeffitients are equal
to zero, and then by (16) in the corresponding local coordinate system the partial
derivatives of tensors $I$ and $J$  coordinates are equal to zero, i.e. the
tensors have constant coordinates in the corresponding local maps.\qed
\edemo
\subsection{Fundamental distributions and integrability}
\definition{Definition 12} Let $(M,I,J)$ be a $\pq$-manifold. On it there
are innerly defined eigendistributions $D^{\l}_I,D^{\mu}_J,D^{\nu}_K$ corresponding to
eigenvalues $\l,\mu,\nu$ of endomor\-phisms $I,J,K=I\circ J$ respectively. We call
the distributions {\it funda\-mental.} The distributions $D^{\l}_I$ are called
{\it the principal fundamental distributions.}
\edefinition
Evidently, $\l$ and $\mu$ are equal to $\pm i$ or $\pm1$, $\nu=\pm i$. The
question naturally arises about the conditions of the distributions being
involutive. It finds its graceful solution in terms of canonical connection.
\par
Let $\nabla$ be canonical connection of $\pq$-structure, $S$ be its torsion
tensor. We introduce in module ${\sf X}(M)$ (as well as in tangent space
$T_p(M),\,p\in M$) a structure of anticommutative algebra with composition
operation $X*Y=S(X,Y);\quad X,Y\in{\sf X}(M)$, and extend it by linearity on
module ${\bf H}_\a\otimes{\sf X}(M)$ (resp., ${\bf H}_\a\otimes T_p(M)$). We call
the algebra {\it adjoint\/} and denote it by $\sf V$ (resp., ${\sf V}_p$).
By ($16_3$),
$$IX*Y=X*IY;\qquad X,Y\in{\sf V}.$$
\procl{Theorem 11} Fundamental distribution of $\pq$-structure is involutive
iff it is a subalgebra of the adjoint algebra. Here the principal fundamental
distribution is involutive iff it is the ideal of the adjoint algebra.
\eprocl
\demo{Proof} Let $F$ be one of structural endomorphisms $I,J$ or $K$, $\l$ be
its eigenvalue, $D^\l_F$ be the corresponding fundamental distribution,
$D^{-\l}_F$ be the complementary fundamental distribution, $\pi^+=\frac12
(\id+\l F^3)$ be projector on $D^\l_F$,\ \ $\pi^-=\frac12(\id-\l F^3)$
be projector on $D^{-\l}_F$. Assume that $D^\l_F$ is involutive, i.e.
$\forall X,Y\in D^\l_F\implies [X,Y]\in D^\l_F$. This condition can be written
in the form $\pi^-[\pi^+X,\pi^+Y]=0;\quad X,Y\in{\sf X}(M)$. In view of
definition of endomorphisms $\pi^+$ and $\pi^-$ it can be rewritten in the form
\begin{equation}
\a N_F(X,Y)-\l F\circ N_F(X,Y)=0;\qquad X,Y\in{\sf X}(M),
\end{equation}
where $N_F(X,Y)=\a[X,Y]+[FX,FY]-F[FX,Y]-F[X,FY]$ is Nijenhuis tensor of the
endomorphism $F$. Evidently, equality (22) can be written in the form
$$\l F^3\circ N_F(X,Y)=N_F(X,Y),\quad\mbox{i.e.}\quad (\pi^-)\circ N_F(X,Y)=0;
\quad X,Y\in{\sf X}(M).$$
In view of $F\circ N_F(X,Y)+N_F(FX,Y)=0,$ that can be immediately checked we
can rewrite (22) also in the form
\begin{equation}
N_F(\pi^+X,Y)=0;\qquad X,Y\in{\sf X}(M).
\end{equation}
On the other hand, by equality $[X,Y]=\n XY-\n YX-S(X,Y)$ and (16) it is easy
to compute that
\begin{eqnarray}
N_F(X,Y)=-\a S(X,Y)-S(FX,FY)+F\circ S(FX,Y)+F\circ S(X,FY)\nonumber\\
=-F(F\circ S(X,Y)-S(FX,Y))+(F\circ S(X,FY)-S(FX,FY)).
\end{eqnarray}
We introduce tensor $U(X,Y)=F\circ S(X,Y)-S(FX,Y)$. Evidently,
$F\circ U(X,Y)+U(FX,Y)=0;\quad X,Y\in{\sf X}(M)$. Using the tensor we can
rewrite (24) as follow $N_F(X,Y)=-F\circ U(X,Y)+U(X,FY)$. Then (23) will be
written in the form $-F\circ U(\pi^+X,Y)+U(\pi^+X,FY)=0$, and thus,
$U(\pi^+X,FY)=F\circ U(\pi^+X,Y)=-U(F\circ\pi^+X,Y)=-U(\l\pi^+X,Y)=U(\pi^+X,
-\l Y)$, i.e. $U(\pi^+X,(F+\l\id)Y)=0$. Note that $F^4=\id$. Then we can rewrite
the equality in the form $U(\pi^+X,(\id+\l F^3)Y)=0,$ that is $U(\pi^+X,\pi^+Y)=0$.
Thus, distribution $D^\l_F$ is involutive iff $U(\pi ^+X,\pi^+Y)=0;\quad
X,Y\in{\sf X}(M)$. But it means that $F\circ S(\pi^+X,\pi^+Y)-S(F\circ\pi^+X,
\pi^+Y)=0$, i.e. $F\circ S(\pi^+X,\pi^+Y)=\l S(\pi^+X,\pi^+Y)$, and therefore,
$S(\pi^+X,\pi^+Y)\in D^\l_F$, i.e. $D^\l_F\subset{\sf V}$ is a subalgebra.
Inversely, if $D^\l_F\subset{\sf V}$ is a subalgebra, then $F\circ S(\pi^+X,
\pi^+Y)=\l S(\pi^+X,\pi^+Y)=S(\l\pi^+X,\pi^+Y)=S(F\circ \pi^+X,\pi^+Y)$, i.e.
$U(\pi^+X,\pi^+Y)=0;\quad X,Y\in{\sf X}(M)$ that is equivalent to distribution
$D^\l_F$ being involutive. If, in particular, $D^\l_F$ is the principal
fundamental distribution, then $\l S(\pi^+X,\pi^-Y)=S(\l\pi^+X,\pi^-Y)=S(F
\circ\pi^+X,\pi^-Y)=S(\pi^+X,F\circ\pi^-Y)=S(\pi^+X,-\l\pi^-Y)=-\l S(\pi^+X,
\pi^-Y)$, i.e. $S(\pi^+X,\pi^-Y)=0;\quad X,Y\in{\sf X}(M)$ and then $S(\pi^+X,
Y)=S(\pi^+X,\pi^+Y)\in D^\l_F$, i.e. $D^\l_F\subset{\sf V}$ is an ideal.\qed
\edemo
\definition{Definition 13} A $\pq$-structure with integrable structural
endomorphism $I$ is called {\it semiholonomic.}
\edefinition
Evidently, the $\pq$-structure being semiholonomic is equivalent to both its
principal fundamental distributions being involutive.
\procl{Theorem 12} $\pq$-structure $\{\,I,J\,\}$ is semiholonomic iff
$$I(X*Y)=I(X)*Y=X*I(Y);\qquad X,Y\in{\sf V}.$$
\eprocl
\demo{Proof} Note that integrability of endomorphism $I$ is equivalent to it
Nijenhuis tensor $N_I(X,Y)$ vanishing, that is, the following equality is valid:
$\a S(X,Y)+S(IX,IY)-I\circ S(IX,Y)-I\circ S(X,IY)=0;\quad X,Y\in{\sf X}(M)$.
The equality can be rewritten by $(16_3)$ in the form $\a S(X,Y)=I\circ S(IX,
Y)$, or $I\circ S(X,Y)=S(IX,Y)=S(X,IY);\quad X,Y\in{\sf X}(M)$.
Inversely, validity of these equalities implies the equality $N_I(X,Y)=0;
\quad X,Y\in{\sf X}(M)$, that is, integrability of endomorphism $I$.\qed
\edemo
\procl{Theorem 13} A non-principal fundamental distribution $D^\l_F$
of semiholo\-nomic $\pq$-structure is involutive iff $FX*FY=\l F(X*Y);\quad
X,Y\in{\sf V}$. In particular, any of the structural endomorphisms $J$ or $K$
of semiholonomic $\pq$-structure is integrable iff the adjoint algebra $\sf V$
is Abelian.\eprocl
\demo{Proof} Since the $\pq$-structure is semiholomorphic, both its principal
fundamen\-tal distributions are involutive. Let $D^\l_F$ be involutive non-principal
fundamental distribution. Then $D^\l_F\subset{\sf V}$ is a subalgebra. Note
that ${\sf V}=D^\mu_I\oplus D^{-\mu}_I$ and that $D^\mu_I*D^{-\mu}_I=\{0\}$.
Indeed, if $X\in D^\mu_I,Y\in D^{-\mu}_I$, then $\mu S(X,Y)=S(\mu X,Y)
= S(IX,Y)=S(X,IY)=S(X,-\mu Y)=-\mu S(X,Y)$, i.e. $S(X,Y)=0.$  Now let $X,Y\in D^\mu_I$.
Then $S(X+\l F^3X,Y+\l F^3Y)=S(X,Y)+\l S(F^3X,Y)+\l S(X,F^3Y)+\l^2S(FX,FY)$.
We show that $F(D^\mu_I)\subset D^{-\mu}_I$. Indeed, if $X\in D^\mu_I$, then
$\mu F(X)=F(\mu X)=F\circ I(X)=-I\circ F(X)$, i.e. $F(X)=-\mu F(X)$ and thus,
$F(X)\subset D^{-\mu}_I$. Consequently, $S(X+\l F^3X,Y+\l F^3Y)=S(X,Y)+\l^2S(FX,FY)$.
But $X+\l F^3X=2\pi^+X\in D^\l_F$, and since $D^\l_F\subset{\sf V}$ is subalgebra,
$S(X,Y)+\l^2S(FX,FY)\in D^\l_F$. But since $D^\mu_I$ and $D^{-\mu}_I$ are also
subalgebras, $S(X,Y)\in D^\mu_I,\ S(FX,FY)\in D^{-\mu}_I$, and hence,
$F\circ S(X,Y)\in D^{-\mu}_I,\ F\circ S(FX,FY)\in D^\mu_I$. In view of
${\bf K}_\a\otimes{\sf X}(M)=D^\mu_I\oplus D^{-\mu}_I$ we have:
$F\circ S(X,Y)=\l^3S(FX,FY)$, or $S(FX,FY)=\l F\circ S(X,Y)$. The same formula
is valid for $X,Y\in D^{-\mu}_I$, and since $S(FX,FY)=0=\lambda F\circ S(X,Y);\quad
X\in D^\mu_I,\quad Y\in D^{-\mu}_I$, we have $S(FX,FY)=\l F(X,Y);\quad X,Y\in
{\sf X}(M)$, i.e. $FX*FY=\l F(X*Y);\quad X,Y\in{\sf V}$. Inversely, if this
equality is hold, then $F(X*Y)=\l^3(FX)*(FY)=\l X*Y;\quad X,Y\in D^\l_F$, i.e.
$D^\l_F\subset{\sf V}$ is a subalgebra, and $D^\l_F$ is involutive distribution.
If $F$ is an integrable endomorphism, then both distributions, $D^\l_F$ and
$D^{-\l}_F$, are involutive, and thus, $F(X*Y)=\l^3FX*FY$ and $F(X*Y)
=-\l^3FX*FY$ simultaneously, hence,
$X*Y=0;\quad X,Y\in{\sf V}$, i.e. $\sf V$ is Abelian algebra. Inversely, if
$\sf V$ is Abelian algebra, then all fundamental distributions are involutive
(since they are automatically subalgebras), in particular, distributions
$D^\l_F$ and $D^{-\l}_F$ are involutive, hence, endomorphism $F$ is integrable.
\qed\edemo
\procl{Corollary 1} Giving three-web on a manifold is equivalent to giving
on it a $\pq$-structure $\{\,I,J\,\}$, its adjoint algebra $\sf V$ being
$J$-linear, and endomorphism I being its involutory automorphism.
\eprocl
\demo{Proof} Let a three-web, i.e. triples of $n$-dimensional involutive pair-wise
comple\-mentary distributions, be given on $2n$-dimensional manifold $M$. The
adjoint almost antiquaternionic structure of the three-web [35] is semiholonomic
with respect to endomorphism $J$, distribution $D^1_I$ being involutive. By
Theo\-rems 12 and 13 its adjoint algebra is $J$-linear, and $I$ its involutory
automor\-phism. Inversely, if $\pq$-structure $\{\,I,J\,\}$ has these properties,
then, by automorphism $I$ being involutory, the $\pq$-structure is antiquaternionic,
and, by the same theorems,  distributions $D^1_J$, $D^{-1}_J$ and $D^1_I$
define on $M$ the structure of three-web.\qed
\edemo
\procl{Corollary 2} Integrability of any pair of structural endomorphisms
of $\pq$-structure $\{\,I,J\,\}$ implies integrability of the third structural
endomor\-phism and is equivalent to the adjoint algebra being Abelian.
\eprocl
\demo{Proof} It follows from the condition that at least one of the endomorphisms,
$I$ or $J$, is integrable. Redenoting, if necessary, $I$ and $J$, we get that
endomorphism $I$ is integrable, i.e. $\pq$-structure is semiholonomic. Since
one of the endomorphisms, $J$ or $K$, is also integrable, by Theorem 13 we have
that the adjoint algebra is Abelian, and thus (by the same Theorem) all
structural endomorphisms are integrable.\qed
\edemo
\remark{Remark} The three-web whose Chern connection is torsion-free and thus
the adjoint algebra is Abelian, is called {\it paratactic,} or {\it isoclinic-geodesic\/} [24].
\eremark
\subsection{Isoclinic distributions}
Let $\{\,M,I,J\,\}$ be a $\pq$-manifold, $\sf Q$ be its structural bundle.
\definition{Definition 14} An $\a$-quaternion $q\in\{\,{\sf Q}\,\}$ whose kernel
forms an $r$-dimensional distribution $(r\ge 0)$ is called {\it admissible.}
\edefinition
Note that if $r>0$, then $q$ is an isotropic $\a$-quaternion, i.e. $|q|=0$.
Indeed, in this case $\forall p\in M\,\exists X\in T_p(M)\mid X\ne 0\ \&\ q(X)=0.$
But then $|q|^2X=\o qq(X)=0$, and thus, $|q|=0$ at arbitrary point $p\in M$,
i.e. $|q|=0$.
\par
The examples of admissible $\a$-quaternions with non-zero kernel are projec\-tors
$\pi^+_F$ and $\pi^-_F$ with respect to structural endomorphism $F$; their kernels
are fundamental distributions $D^{-\l}_F$ and $D^\l_F$, respectively. Now we set
up the task of studying the kernel of arbitrary  admissible $\a$-quaternions.
Further on all $\a$-quaternions being regarded are assumed admissible.
\definition{Definition 15} The kernel of an isotropic $\a$-quaternion is called
{\it an isoclinic distribution\/} on $\pq$-manifold.
\edefinition
Let $q$ be an isotropic $\a$-quaternion, $q=a\id+bI+cJ+dK,\,D=\ker q$ is the
corresponding isoclinic distribution, $X\in D$. Note that ${\bf K}_\a\otimes
{\sf X}(M)=D^\l_I\oplus D^{-\l}_I$, and thus, $X=X^++X^-,\quad X^+\in D^\l_I,
\,X^-\in D^{-\l}_I$. Then $q(X)=(a\id+bI+cJ+dK)(X^++X^-)=(a+b\l)X^++(c+d\l)
JX^-+(c-d\l)JX^++(a-b\l)X^-=0$, hence, $(a+b\l)X^++(c+d\l)JX^-=0,\
(a-b\l)X^-+(c-d\l)JX^+=0$, or
$$\left\{\begin{array}{rcl}
(a+b\l)X^++(c+d\l)JX^-\eq0;\\
\a(c-d\l)X^++(a-b\l)JX^-\eq0.
\end{array}\right.$$
The determinant of the given system of equations is equal to $|q|^2=0$, and
thus, vectors $X^+$ and $JX^-$ are collinear. Let $X^+\ne 0$. Then $JX^-=\nu X^+$,
where $\nu=-\frac{a+b\l}{c+d\l}$ or $\nu=-\a\frac{c-d\l}{a-b\l}$. Similarly,
if $X^-\ne 0$, i.e. $X\notin D^\l_I$, then $JX^+=\nu X^-$, where
$\nu=-\a\frac{c+d\l}{a+d\l}$ or $\nu=-\frac{a-b\l}{c-d\l}$. Thus,
$$D=\{\,X+\mu JX\mid X\in D^\l_I\,\},$$
where $\mu\in C^\infty\,(M)$ is fixed function, $D^\l_I$ is one of the principal
fundamental distributions. Distributions of such kind are called {\it isoclinic\/}
in three-web theory [37].
\par
Inversely, consider distribution of the type $D_\mu=\{\,X+\mu JX\mid\mu\in
C^\infty(M),\ X\in D^\l_I\,\}$ on a $\pq$-manifold $M$, and show, that it is
isoclinic. If $\ \mu\ne$ const, we call it {\it\ a nontrivial isoclinic
distribution.} We find all $\a$-quaternions $q$, such that $D_\mu=\ker q$.
Let $q=a\id+bI+cJ+dK$. Then
$$q(X+\mu IX)=0\iff\left\{\begin{array}{rcl}\a\mu(c+d\l)\eq-(a+b\l),\\
\mu(a-b\l)\eq-(c-d\l),
\end{array}\right.$$
hence, $a=-\frac12(\frac1\mu+\a\mu)c+\frac12(\frac1\mu-\a\mu)d,\
b=\frac1{2\l}(\frac1\mu-\a\mu)c-\frac12(\frac1\mu+\a\mu)d$, and thus,
\begin{eqnarray*}
q\eq\{-\frac12(\frac1\mu+\a\mu)c+\frac12(\frac1\mu-\a\mu)d\,\}\id\,+\\
&{}&+\{\,\frac1{2\l}(\frac1\mu-\a\mu)c-\frac12(\frac1\mu+\a\mu)d\,\}\,I+cJ+dK\\
\eq\{-\frac12(\frac1\mu+\a\mu)\id+\frac1{2\l}(\frac1\mu-\a\mu)I+J\,\}\,c\,+\\
&{}&+\{\,\frac12(\frac1\mu-\a\mu)\id-\frac12(\frac1\mu+\a\mu)I+K\,\}\,d\,.
\end{eqnarray*}
Thus, distribution $D_\mu$ defined as a common kernel of $\a$-quaternions
forming a 2-dimensional submodule of structural bundle section module. Thus,
we have proved
\procl{Theorem 14} There exists a natural one-to-one correspondence between
isoclinic distributions of a $\pq$-structure and 2-dimensional distributions
of isotropic $\a$-quaternions of the mentioned type on a $\pq$-manifold.\qed
\eprocl
It follows from the above that a unification of isoclinic distributions on a
$\pq$-manifold gives the foliation of cones over the manifold modelled by
isotropic cones of corresponding $\a$-quaternions. The elements of this foliation
will be called {\it Segre cones\/} (as in three-web theory), and its points are
called {\it isoclinic vectors}. Isoclinic subspaces at a given point of a
$\pq$-manifold are $\frac12\dim M$-dimensional rulings of Segre cones at the point.
\definition{Definition 16} A $\pq$-manifold will be called {\it isoclinic} if
it admits an involutive isoclinic distribution different from fundamental distributions.
\edefinition
Let $\{\,M,I,J\,\}$ be a semiholonomic isoclinic $\pq$-manifold, $q\in{\sf Q}$
be the isotropic $\a$-quaternion defining involutive isoclinic distribution
$D_\mu$. Since $D_\mu$ is involutive, $\forall X,Y\in D_\mu\implies Z=[X,Y]\in D_\mu$.
By the proved above, $X=a+\mu Ja$, $Y=b+\mu Jb$, $Z=c+\mu Jc$, where $a,b,c\in D^\l_I$.
Thus, $[X,Y]=[a+\mu Ja,b+\mu Jb]=[a,b]+\mu([a,Jb]+ [Ja,b])+\mu^2[Ja,Jb]+a(\mu)Jb
+\mu Ja(\mu)Jb-b(\mu)Ja-\mu Jb(\mu)Ja=c+\mu Jc$.
Note that $N_J(a,b)=\a[a,b]+[Ja,Jb]-J[Ja,b]-J[a,Jb]$ and then
$[a,Jb]+[Ja,b]=\a J(-N_J(a,b)+\a[a,b]+[Ja,Jb])=\a J(\a S(a,b)+S(Ja,Jb)
+\a[a,b]+[Ja,Jb])=J\circ S(a,b)+\a J\circ S(Ja,Jb)+J[a,b]+\a J[Ja,Jb]$.
Thus, $[a,b]+\mu J\circ S(a,b)+\a\mu J\circ S(Ja,Jb)+\mu J[a,b]+\a\mu J[Ja,Jb]
+\mu^2[Ja,Jb]+a(\mu)Jb+\mu Ja(\mu)Jb-b(\mu)Ja-\mu Jb(\mu)Ja=c+Jc$.
In view of ${\bf K}_\a\otimes{\sf X}(M)=D^\l_I\oplus D^{-\l}_I$, and distributions
$D^\l_I$ and $D^{-\l}_I$ being involutive and, thus, being subalgebras of the
adjoint algebra, we get:
$$\begin{array}{rl}
&1)\ [a,b]+\a\mu JS(Ja,Jb)+\a\mu J[Ja,Jb]=c,\\
&2)\ \mu JS(a,b)+\mu J[a,b]+\mu^2[Ja,Jb]+a(\mu)Jb\,+\\
&\quad+\mu Ja(\mu)Jb-b(\mu )Ja-\mu Jb(\mu)Ja=\mu Jc,
\end{array}$$
and thus,
\begin{eqnarray*}
1)\ c\eq[a,b]+\a\mu JS(Ja,Jb)+\a\mu J[Ja,Jb];\\
2)\ c\eq[a,b]+S(a,b)+\a\mu J[Ja,Jb]+\frac1\mu
a(\mu)b+Ja(\mu)b-\frac1\mu b(\mu)a-Jb(\mu)a.
\end{eqnarray*}
Substract elementwise one equality from the other:
$S(a,b)-\a\mu J\circ S(Ja,Jb)=(\frac1\mu b(\mu)+Jb(\mu))a-(\frac1\mu a(\mu)+Ja(\mu))b$,
or $J\circ S(a,b)-\mu S(Ja,Jb)=\omega(b)Ja-\omega(a)Jb$, where $\omega(c)
=\frac1\mu(c+\mu Jc)(\mu)=\frac1\mu d\mu\,(c+\mu Jc)=d\,(\ln|\mu|)\circ\id+\mu J)c$.
Evidently, the inverse is also true: this condition implies
that isoclinic distribution $D_\mu$ is involutive. Thus, we have proved
\procl{Theorem 15} A semiholonomic $\pq$-manifold $\{\,M,I,J\,\}$ is
isoclinic iff
$$\exists\mu\in C^\infty\,(M);\mu\ne\pm 1\ \&\ J(X*Y)-\mu(JX*JY)=\omega(Y)JX-\omega(X)JY;$$
where $X,Y\in D^\l_I,\ \omega=d\,(\ln|\mu|)(\id+\mu J)$. In this case
$D=\{\,X+\mu JX\mid X\in D^\l_I\,\}$ is an involutive isoclinic distribution on M.\qed
\eprocl
\procl{Corollary 1} A semiholonomic $\pq$-manifold $\{\,M,I,J\,\}$ having
not less than three involutive fundamental distributions is isoclinic iff
$$X*Y=\frac1{1-\mu}\{\,\omega(Y)X-\omega(X)Y\,\};\qquad X,Y\in D^\l_I.\qeds$$
\eprocl
The above result generalizes the known criterion of M.A.Akivis of three-web
being isoclinic [37].
\procl{Corollary 2} A semiholonomic $\pq$-manifold having more than three
involutive fundamental distributions is isoclinic.\qed
\eprocl
\definition{Definition 17} A $\pq$-manifold is called {\it isoclinic-geodesic\/}
if it admits a non-trivial totally geodesic (in canonical connection) isoclinic
distribution.
\edefinition
Let $\{\,M,I,J\,\}$ be a $\pq$-manifold, $q\in\{\,{\sf Q}\,\}$ be an isotropic
$\a$-quaternion defining isoclinic distribution $D=\{\,X+\mu JX\mid X\in
D^\l_I\,\}$. Let $X,Y\in D,\ X=a+\mu Ja,\ Y=b+\mu Jb,\ a,b\in D^\l_I$. Then
$\n XY=\n ab+\mu\nabla_{Ja}\,b+\mu\n a(Jb)+\mu^2\nabla_{Ja}\,(Jb)+a(\mu)Jb+\mu Ja(\mu)Jb
=\n ab+\mu J\n ab+\mu(\nabla_{Ja}\,b+\mu J\nabla_{Ja}\,b)+a(\mu)Jb+\mu Ja(\mu)Jb$.
It is clear that $D$ is totally geodesic, i.e. $\n XY\in D\quad(X,Y\in D)$
iff $a(\mu)+\mu Ja(\mu)=0$, i.e., in the notations of the above theorem,
$\omega=0$. We have proved
\procl{Theorem 16} A $\pq$-manifold is isoclinic-geodesic iff $\omega=0$
where
$$\omega=d\,(\ln|\mu|)\circ(\id+\mu J);\quad\mu\in C^\infty\,(M);\ \mu\ne const.$$
In this case $D=\{\,X+\mu JX\mid X\in D^\l_I\,\}$ is a totally geodesic
isoclinic distribution on this manifold.\qed
\eprocl
Now let $\{\,M,I,J\,\}$ be a semiholonomic $\pq$-manifold, $D$ be an involutive
isoclinic distribution on $M$, $X,Y\in D$, $X=a+\mu Ja$, $Y=b+\mu Jb$; $a,b\in D^\l_I$.
Then $X*Y=a*b+\mu(Ja*b+a*Jb)+\mu^2(Ja*Jb)$. In view of Theorem 11, $D^\l_I$
and $D^{-l}_I$ are ideals of the adjoint algebra, and since $J(D^\l_I)\subset D^{-\l}_I$,
$Ja*b=a*Jb=0$. Thus, $X*Y\in D\iff a*b+\mu^2(Ja*Jb)=c+\mu Jc;
\ (c\in D^\l_I)$, hence, $c=a*b,Jc=\mu(Ja*Jb)$, and consequently, $J(a*b)=\mu(Ja*Jb);
\ (a,b\in D^\l_I)$. Evidently, the inverse is also true, i.e. there holds
\procl{Theorem 17} An involutive isoclinic distribution $D=\{\,X+\mu JX\mid X\in D^\l_I\,\}$
of a semiholonomic $\pq$-manifold $M$ is a subalgebra of the adjoint algebra
iff $$J(X*Y)=\mu(JX*JY);\qquad X,Y\in D^\l_I.\qeds$$
\eprocl
Further, by Theorems 15 and 16 we get
\procl{Theorem 18} A semiholonomic isoclinic $\pq$-manifold $M$ is
isoclinic-geodesic iff $J(X*Y)=\mu(JX*JY);\quad X,Y\in D^\l_I;\ \mu\in C^\infty\,(M);
\mu\ne const$. In this case $D=\{\,X+\mu JX\mid X\in D^\l_I\,\}$ is a totally
geodesic isoclinic distribution on $M$ being a subalgebra of the adjoint algebra.\qed
\eprocl
\procl{Corollary} A semiholonomic isoclinic $\pq$-manifold having not less
than three involutive fundamental distributions is isoclinic-geodesic iff its
adjoint algebra is Abelian.
\eprocl
\demo{Proof} It immediately follows from Theorem 16 and Corollary 1 of Theorem 15.\qed
\edemo
The result shows that the notion of an isoclinic-geodesic $\pq$-structure is
the generalization of the notion of an isoclinic-geodesic three-web, mentioned
above and, moreover, it explains the geometric sense of the notion.
\subsection{Hermitian and Einsteinian metrics generated by anti\-quaternionic
structure on pseudo-Riemannian mani\-fold}
Let $(M,(\cdot,\cdot))$ be a pseudo-Rie\-mannian manifold, $\{\,I,J\,\}$ be a
$\pq$-structure on $M$. We define the pseudo-Riemannian metric $g$ on $M$ by
the formula $$g(X,Y)=(X,Y)+(IX,IY)+(JX,JY)+(KX,KY);\qquad X,Y\in{\sf X}(M).$$ Evidently,
\begin{eqnarray}
g(IX,IY)=g(X,Y);\qquad g(IX,Y)\eq\a g(X,IY);\nonumber\\
g(JX,JY)=g(X,Y);\qquad g(JX,Y)\eq\a g(X,JY);\\
g(KX,KY)=g(X,Y);\qquad g(KX,Y)\eq-g(X,KY)\nonumber.
\end{eqnarray}
Let $\a=+1$, i.e. the $\pq$-structure under consideration is almost antiquaterni\-onic.
Then two neutral metrics $g_1(X,Y)=g(X,IY)$ and $g_2(X,Y)=g(X,JY)$ are defined
on $M$ in addition. Indeed,
\begin{eqnarray*}
&g_1(JX,JY)=g(JX,IJY)=-g(JX,JXY)=-g(X,IY)=-g_1(X,Y);\\
&g_2(IX,IY)=g(IX,JIY)=-g(IX,IJY)=-g(X,JY)=-g_2(X,Y).
\end{eqnarray*}
Thus, whole metrics sheaf
\begin{equation}
\la X,Y\ra=g(X,Y)+\l g(X,IY)+\mu g(X,JY)
\end{equation}
is associated with $g$. It is easy to see that the metrics are non-degenerate
iff $\l^2+\mu^2\ne1$. Indeed, let $\la X,Y\ra=0;\quad X,Y\in{\sf X}(M)$.
Then by non-degeneracy of $g$ we have:
\begin{equation}
X+\l IX+\mu JX=0.
\end{equation}
Acting on both parts of the equality at first by endomorphism $I$ and then by
$J$  we have respectively:
$$\l X+IX+\mu KX=0;\qquad\mu X+JX-\l KX=0.$$
We multiply both parts of the first equality on $\l$, and of the second one by
$\mu$ and add the received equalities term by term:
$$(\l^2+\mu^2)X+\l IX+\mu JX=0.$$
In view of  (27) we see that either $\l^2+\mu^2=1$, or $X=0$. We find all purely
imaginary antiquaternions $J$ forming almost Hermitian structure with a metric
of such kind. Let $J=\b I+\c J+\d K$. By direct calculations we find that
$\{\,J,\la\cdot,\cdot\ra\,\}$ is almost Hermitian structure on $M$, i.e.
$\la JX,JY\ra=\la X,Y\ra;\quad X,Y\in{\sf X}(M)$ iff the following equation
system is fulfilled:
$$\left\{\begin{array}{l}
\b^2+\c^2+\d^2+2\mu\b\c=1;\\
-2\c\d+2\mu\b\c+\l(\b^2-\c^2-\d^2)=\l;\\
2\b\d+\mu(-\b^2+\c^2-\d^2)+2\l\b\c=\mu;\\
-\b^2-\c^2+\d^2=1.
\end{array}\right.$$
According to the last system equation, $\b^2+\c^2=\d^2-1$. Therefore, the system
can be rewritten in the form:
$$\left\{\begin{array}{l}
\d\varphi=1;\\\c\varphi=-\l;\\\b\varphi=\mu;\\\b^2+\c^2=\d^2-1,
\end{array}\right.$$
where $\varphi=\d-\mu\b+\l\c$. We have from the first triple equations:
$$\varphi^2(-\b^2-\c^2+\d^2)=1-\l^2-\mu^2.$$
Note that in view of the first system equation $\d\ne0$ and hence,
$\frac\c\d=-\l,\ \frac\b\d=\mu$, i.e. $\b=\mu\d,\ \c=-\l\d$, and in view of
the last system equation,
$$
\b=\pm\frac\mu{\sqrt{1-\l^2-\mu^2}};\quad\c=\pm\frac\l{\sqrt{1-\l^2-\mu^2}};
\quad\d=\pm\frac1{\sqrt{1-\l^2-\mu^2}}.
$$
Thus we have proved
\procl{Theorem 19} Every pseudo-Riemannian metric $(\cdot,\cdot)$ on $\aq$-manifold
$M$ is associated with two-parameter family of almost Hermitian structures
$\{\,{\cal J},\la\cdot,\cdot\ra\,\}$ on $M$ ;
$${\cal J}=\pm\frac{\mu I-\l J+K}{\sqrt{1-\l^2-\mu^2}};\qquad
\la X,Y\ra=g(X,Y)+\l g(X,IY)+\mu g(X,JY);$$
$$\l,\mu\in{\bf R},\quad\l^2+\mu^2<1;$$
where $g(X,Y)=(X,Y)+(IX,IY)+(JX,JY)+(KX,KY);\ X,Y\in{\sf X}(M)$.\qed
\eprocl
\remark{Remark} In the case $\l^2+\mu^2>1$ every pseudo-Riemannian metric is
associated with two-parameter family of hyperbolic almost Hermitian structures
$\{\,{\cal J},\la\cdot,\cdot\ra\,\}$, where
$${\cal J}=\pm\frac{\mu I-\l J+K}{\sqrt{\l^2+\mu^2-1}},\qquad
\la X,Y\ra=g(X,Y)+\l g(X,IY)+\mu g(X,JY).$$
Indeed, in the case ${\cal J}^2=\id$, i.e. $\|{\cal J}\|^2=-1$, and therefore,
$-\b^2-\c^2+\d^2=-1$, hence, $\d=\pm\frac1{\sqrt{\l^2+\mu^2-1}}$.
\eremark
Now let $M$ be parallelizable locally homogeneous naturally reductive mani\-fold,
$(\cdot,\cdot)$ be invariant pseudo-Riemannian metric on $M$, $G$ be fundamental
group of local isometries of $M$, $H\subset G$ be isotropy subgroup, $\sf g$
and $\sf h$ be Lie algebras of Lie groups $G$ and $H$, respectively. Then
${\sf g}={\sf h}\oplus{\sf m}$, $\ad({\sf h}){\sf m}\subset{\sf m}$. As usual,
we shall identify a tangent space $T_p(M);\quad p\in M$, and subspace
${\sf m}\subset{\sf g}$ in canonical way. Natural reductivity means that
$([X,Y]_{\sf m},Z)=(X,[Y,Z]_{\sf m});\quad X,Y,Z\in{\sf m}$. Recall [35] that
manifold $M\times M$ carries naturally defined $\aq$-structure $\{\,I,J\,\}$.
The distributions $D^1_I$ and $D^{-1}_I$ of the structure are generated by
subspaces ${\sf m}_1$ and ${\sf m}_2$ respectively, where ${\sf m}_i$ is the
$i$-th term of direct sum ${\sf m}\oplus{\sf m}$ corresponding to decomposition
$T_{(p,q)}(M\times M)=T_p(M)\oplus T_q(M)$. Since distributions $D^1_I$
and $D^{-1}_I$ are involutive the $\aq$-structure is semiholonomic. Moreover,
let us denote the diagonal of manifold $M\times M$ by $\Delta$. Then the distribution
$D^1_J$ is generated by subspace ${\sf m}_0$ corresponding to $T_p(\Delta);
\quad p\in\Delta$, hence, it is involutive also. Therefore, for the stated
$\aq$-structure we have:
\begin{equation}
I(X)*Y=X*I(Y)=I(X*Y);\qquad J(X*Y) = J(X)*J(Y);
\end{equation}
where $X,Y\in{\sf X}(M\times M)$. Let us calculate a torsion tensor $S$ of the
canonical connection $\tilde\nabla$ for the $\aq$-structure. Note that
$[X,Y]=0;\quad X\in{\sf m}_1,Y\in{\sf m}_2$, and in view of (19) we have:
$\tilde\n XY=0;\quad X,Y\in{\sf m}\oplus{\sf m}$. Therefore, $X*Y=S(X,Y)
=\tilde\n XY-\tilde\n YX-[X,Y]=-[X,Y];\quad X,Y\in{\sf m}\oplus{\sf m}$. Here
and later on by $[X,Y]$ we mean of vectors $X$ and $Y$  commutator restriction
on ${\sf m}\oplus{\sf m}\subset{\sf g}\oplus{\sf g}$. Consequently, relations
(28) take the form:
\begin{equation}
[I(X),Y]=[X,I(Y)]=I([X,Y]);\qquad J([X,Y])=[J(X),J(Y)];
\end{equation}
where $X,Y\in{\sf m}\oplus{\sf m}$.
\par
We calculate Riemannian connection $\nabla$ of the metric $\o g=\la\cdot,\cdot\ra$
generated by direct product metric of manifold $M\times M$ as mentioned previously. We have:
$$2\la\n XY,Z\ra=\la[X,Y],Z\ra+\la[Z,X],Y\ra-\la[Y,Z],X\ra.$$
After direct calculations with respect to natural reductivity of manifold we get:
\begin{equation}
Z+\l IZ+\mu JZ=A,
\end{equation}
where $Z=2\,\n XY$, $A=[X,Y]+\mu J[X,Y]+\mu[X,JY]-\mu[JX,Y]+\l I[X,Y]$.
Acting on (30) by endomorphisms $I$, $J$ and $K$ we receive correspondingly:
\begin{eqnarray}
&1)\ \l Z+IZ+\mu KZ=IA;\nonumber\\
&2)\ \mu Z+JZ-\l KZ=JA;\\
&3)\ KZ-\l JZ+\mu IZ=KA\nonumber.
\end{eqnarray}
Multiply the third equation in (31) first by $\mu$ and then by $\l$:
\begin{eqnarray}
&1)\ \mu KZ-\mu\l JZ+\mu^2IZ=\mu KA;\nonumber\\
&2)\ \l KZ-\l^2JZ+\l\mu IZ=\l KA.
\end{eqnarray}
Subtract $(32_1)$ from $(31_1)$ termwise and add $(32_2)$ and $(31_2)$ termwise:
\begin{eqnarray}
&1)\ \l Z+(1-\mu^2)IZ+\l\mu JZ=(I-\mu K)A;\nonumber\\
&2)\ \mu Z+\l\mu IZ+(1-\l^2)JZ=(J+\l K)A.
\end{eqnarray}
Now multiplying (30) first by $\l$ and then by $\mu$ we receive, correspondingly:
\begin{eqnarray}
&1)\ \l Z+\l^2IZ+\l\mu JZ=\l A;\nonumber\\
&2)\ \mu Z+\l\mu IZ+\mu^2JZ=\mu A;
\end{eqnarray}
Subtract $(34_1)$ from $(33_1)$ and $(34_2)$ from $(33_2)$ termwise:
$$Z=\frac{\id-\l I-\mu J}{1-\l^2-\mu^2}A,$$
therefore,
$$
\n XY=\frac12\frac{\id-\l I-\mu J}{1-\l^2-\mu^2}\{\,[X,Y]
+\mu J[X,Y]+\mu[X,JY]-\mu[JX,Y]+\l[X,Y]\,\}.
$$
Simplifying the equation with regard to (29) we get finally:
\begin{eqnarray}
\n XY\eq\frac12[X,Y]+\frac{\mu^2+\mu}{2(1-\l^2-\mu^2)}([X,JY]-[JX,Y])+\nonumber\\
&{}&+\frac{\l\mu}{2(1-\l^2-\mu^2)}([KX,Y]-[X,KY]).
\end{eqnarray}
With respect to (35) we calculate covariant differential of tensors $I,J,K$
and $\cal J$ in Riemannian connection: $\n X(I)Y=\n X(IY)-I(\n XY)$. After
direct calculations with regard to (29) we get:
$$\n X(I)Y=-\frac{\mu^2+\mu}{1-\l^2-\mu^2}[X,KY]+\frac{\l\mu}{1-\l^2-\mu^2}[X,JY].$$
Similarly,
\begin{eqnarray*}
\n X(J)Y\eq\frac12\{\,[X,JY]-J[X,Y]\,\}+\\
&{}&+\frac{\mu^2+\mu}{2(1-\l^2-\mu^2)}\{\,[X,Y]-J[X,Y]-[JX,Y]+[X,JY]\,\}+\\
&{}&+\frac{\l\mu}{2(1-\l^2-\mu^2)}\{-I[X,Y]+K[X,Y]-[KX,Y]+[X,KY]\,\};\\
\n X(K)Y\eq\frac12\{\,[X,KY]-K[X,Y]\,\}+\\
&{}&+\frac{\mu^2+\mu}{2(1-\l^2-\mu^2)}\{-I[X,Y]-K[X,Y]-[KX,Y]+[X,KY]\,\}+\\
&{}&+\frac{\l\mu}{2(1-\l^2-\mu^2)}\{\,[X,Y]+J[X,Y]-[JX,Y]+[X,JY]\,\}.
\end{eqnarray*}
Taking into account the definition of tensor $\cal J$ we receive from here:
\begin{eqnarray*}
\n X({\cal J})Y\eq\frac1{2(1-\l^2-\mu^2)}\{\,\l\mu^2[X,Y]+(\mu\l^2-\mu^2-\mu)I[X,Y]+\\
&{}&+(\l-\l^3+2\l\mu)J[X,Y]+(\l^2-\l^2\mu-\mu-1)K[X,Y]+\\
&{}&+(\l^3+2\l\mu^2-\l)[X,JY]+\l\mu^2[JX,Y]+(1+\mu-2\mu^3-\\
&{}&-2\mu^2-\mu\l^2-\l^2)[X,KY]+(\l^2\mu-\mu^2-\mu)[KX,Y]\,\}.
\end{eqnarray*}
In particular, the structure $\{\,{\cal J},\la\cdot,\cdot\ra\,\}$ is nearly Kaehlerian,
i.e. $\n X({\cal J})X=0;\quad X\in{\sf X}(M)$ iff
$$(\l^3+\l\mu^2-\l)[X,JX]+(1+2\mu-2\mu^3-\mu^2-2\l^2\mu-\l^2)[X,KX]=0,$$
i.e.
$$\left\{\begin{array}{l}
1)\ \l(\l^2+\mu^2-1)=0;\\
2)\ (1+2\mu)(\l^2+\mu^2-1)=0,
\end{array}\right.$$
and since $\l^2+\mu^2\ne1$, we have: $\l=0,\ \mu=-\frac12$.
We get the following result:
\procl{Theorem 20} Almost Hermitian structure $\{\,{\cal J},\la\cdot,\cdot\ra\,\}$
on manifold $M\times M$, des\-cribed above is nearly Kaehlerian iff $\l=0,\mu=-\frac12$,
i.e.
$${\cal J}=\pm\frac{I-2K}{\sqrt3};\qquad\la X,Y\ra=g(X,Y)-\frac12g(X,JY).\qeds$$
\eprocl
Now we calculate the compositional tensor of adjoint $Q$-algebra [28] of almost
Hermitian structure $\{\,{\cal J},\la\cdot,\cdot\ra\,\}$ on the manifold
$M\times M$ described above. After direct calculations we get:
$\n X({\cal J}){\cal J}Y+\nabla_{{\cal J}X}({\cal J})Y=(2b\mu-e\l-f\l-h-r)
[X,Y]+(2a\mu-h\l-r\l-e-f)I[X,Y]-(2d\mu-e\l+f\l-h-r)J[X,Y]-(2c\mu+h\l+r\l-f
-e)K[X,Y]-(c\l+a\l-b+d)[JX,Y]-(a\l+c\l-b+d)[X,JY]-(d\l+b\l-a+c)[KX,Y]
-(b\l+d\l-a+c)[X,KY]$, where
\begin{eqnarray*}
a=-\frac{\l\mu^2}{2\sqrt{(1-\l^2-\mu^2)^3}};\
b=\frac{\l^2\mu-\mu^2-\mu}{2\sqrt{(1-\l^2-\mu^2)^3}};\
c=\frac{\l-\l^3+2\l\mu}{2\sqrt{(1-\l^2-\mu^2)^3}};\\
d=\frac{\l^2-\l^2\mu-\mu-1}{2\sqrt{(1-\l^2-\mu^2)^3}};\
e=\frac{\l^3+2\l\mu^2-\l}{2\sqrt{(1-\l^2-\mu^2)^3}};\
f=\frac{\l\mu^2}{2\sqrt{(1-\l^2-\mu^2)^3}};\\
h=\frac{1+\mu-2\mu^3-2\mu^2-\mu\l^2-\l^2}{2\sqrt{(1-\l^2-\mu^2)^3}};\
r=\frac{\mu\l^2-\mu^2-\mu}{2\sqrt{(1-\l^2-\mu^2)^3}}.\end{eqnarray*}
In particular, the tensor is skew-symmetric by covariant arguments, hence,
almost Hermitian structure $\{\,{\cal J},\la\cdot,\cdot\ra\,\}$ belongs to class
$G_1$ in Gray-Hervella classification [28]. Note also, that
$2b\mu-e\l-f\l-h-r=(\l^2+\mu^2-1)(\l^2-1)$. Since $\l^2\le\l^2+\mu^2<1$,
the coefficient at $[X,Y]$ is nonzero, hence, almost Hermitian structure
$\{\,{\cal J},\la\cdot,\cdot\ra\,\}$ is nonintegrable, if $M$ differs from locally
symmetric space. Besides we get that $\{\,{\cal J},\la\cdot,\cdot\ra\,\}$ is
quasi-Kaehlerian structure, i.e. $\nabla_{{\cal J}X}\,({\cal J}){\cal J}Y
+\n X({\cal J})Y=0$,
iff it is nearly Kaehlerian, i.e. $\l=0$, $\mu=-\frac12$. Formulate the results obtained:
\procl{Theorem 21} Let $\{\,M,(\cdot,\cdot)\,\}$ be a parallelizable locally homogeneous
natu\-rally reductive Riemannian manifold. Then the canonical $\aq$-structure
inherently generates two-parameter family of almost Hermitian structures on
manifold $M\times M$. The structures belong to class $G_1$. If $M$ differs from
locally symmetric space the structures are nonintegrable. The family contains
the only (up to sign) quasi-Kaehlerian structure $\{\,{\cal J}_0,\la\cdot,\cdot\ra_0\,\}$
at $\l=0,\,\mu=-\frac12$, which is a nearly Kaehlerian structure.\qed
\eprocl
Theorem 21 gives a wide spectrum of examples of nonintegrable almost Hermitian
$G_1$-structures, as well as examples of proper (i.e. different from Kaehlerian)
nearly Kaehlerian structures. Note that up to now only a limited number of
examples of the mentioned structures is known.
\par
Returning to considering metrics $\la X,Y\ra=g(X,Y)+\l g(X,IY)+\mu g(X,JY)$
on manifold $M\times M$ we calculate Riemann-Christoffel tensor
$$R(X,Y)Z=\nabla_X\,\n YZ-\nabla_Y\,\n XZ-\nabla_{[X,Y]}\,Z$$ of the metrics.
Direct though cumbersome calculations taking into account (35) give the following result:
\begin{eqnarray}
R(X,Y)Z\eq\frac14[[X,Y],Z]+\frac{\mu^2+\mu}{4(1-\l^2-\mu^2)}\{\,[X,[JY,JZ]]
-[Y,[JX,JZ]]\,\}-\nonumber\\
&{}&-\frac{\mu^2+\mu}{4(1-\l^2-\mu^2)}\{\,[JX,[Y,Z]]-[JY,[X,Z]]\,\}+\nonumber\\
&{}&+\frac{\l\mu}{4(1-\l^2-\mu^2)}\{\,[KX,[Y,Z]]-[KY,[X,Z]]\,\}-\nonumber\\
&{}&-\frac{\l\mu}{4(1-\l^2-\mu^2)}\{\,[X,[KY,JZ]]-[Y,[KX,JZ]]\,\}+\nonumber\\
&{}&+\frac{\mu^2+\mu}{4(1-\l^2-\mu^2)}\{\,[X,[Y,JZ]]-[Y,[X,JZ]]\,\}+\nonumber\\
&{}&+\frac{(\mu^2+\mu)^2}{4(1-\l^2-\mu^2)^2}\{\,[X,[JY,Z]]-[Y,[JX,Z]]\,\}-\nonumber\\
&{}&-\frac{(\mu^2+\mu)^2}{4(1-\l^2-\mu^2)^2}\{\,[JX,[Y,JZ]]-[JY,[X,JZ]]\,\}+\nonumber\\
&{}&+\frac{(\mu^2+\mu)\mu\l}{4(1-\l^2-\mu^2)^2}\{\,[KX,[Y,JZ]]-[KY,[X,JZ]]\,\}-\nonumber\\
&{}&-\frac{\mu^2+\mu}{4(1-\l^2-\mu^2)}\{\,[X,[JY,Z]]-[Y,[JX,Z]]\,\}-\nonumber\\
&{}&-\frac{(\mu^2+\mu)^2}{4(1-\l^2-\mu^2)^2}\{\,[X,[Y,JZ]]-[Y,[X,JZ]]\,\}+\nonumber\\
&{}&+\frac{(\mu^2+\mu)^2}{4(1-\l^2-\mu^2)^2}\{\,[JX,[JY,Z]]-[JY,[JX,Z]]\,\}-\\
&{}&-\frac{(\mu^2+\mu)\mu\l}{4(1-\l^2-\mu^2)^2}\{\,[KX,[JY,Z]]-[KY,[JX,Z]]\,\}+\nonumber\\
&{}&+\frac{\l\mu}{4(1-\l^2-\mu^2)}\{\,[X,[KY,Z]]-[Y,[KX,Z]]\,\}+\nonumber\\
&{}&+\frac{(\l\mu)^2}{4(1-\l^2-\mu^2)^2}\{\,[KX,[KY,Z]]-[KY,[KX,Z]]\,\}-\nonumber\\
&{}&-\frac{\mu^2+\mu}{2(1-\l^2-\mu^2)}\{\,[[X,Y],JZ]-[[JX,JY],Z]\,\}+\nonumber\\
&{}&+\frac{\l\mu}{2(1-\l^2-\mu^2)}\{\,[[X,Y],KZ]-[[KX,JY],Z]\,\}-\nonumber\\
&{}&-\frac{(\mu^2+\mu)\mu\l}{4(1-\l^2-\mu^2)}\{\,[JX,[KY,Z]]-[JY,[KX,Z]]\,\}+\nonumber\\
&{}&+\frac{(\l\mu)^2}{4(1-\l^2-\mu^2)^2}\{\,[X,[Y,JZ]]-[Y,[X,JZ]]\,\}-\nonumber\\
&{}&-\frac{\l\mu}{4(1-\l^2-\mu^2)}\{\,[X,[Y,KZ]]-[Y,[X,KZ]]\,\}+\nonumber\\
&{}&+\frac{(\mu^2+\mu)\mu\l}{4(1-\l^2-\mu^2)}\{\,[JX,[Y,KZ]]-[JY,[X,KZ]]\,\}-\nonumber\\
&{}&-\frac{(\l\mu)^2}{4(1-\l^2-\mu^2)^2}\{\,[KX,[Y,KZ]]-[KY,[X,KZ]]\,\}-\nonumber\\
&{}&-\frac{(\l\mu)^2}{4(1-\l^2-\mu^2)^2}\{\,[X,[JY,Z]]-[Y,[JX,Z]]\,\}.\nonumber
\end{eqnarray}
With regard to the identity we shall calculate Ricci tensor $r(X)=g^{ij}\,R(X,e_i)e_j$
of the manifold. We choose a base $\{\,e_a,e_{\hat a}\,\}$ of a space
$T_p(M\times M),\quad p\in M$, where $e_a\in(D^1_I)_p$, $e_{\hat a}\in(D^{-1}_I)_p$,
$e_{\hat a}=J(e_a),\quad a=1,\ldots,n=\dim M,\ \{\,e_1,\ldots$ $\ldots,e_{n}\,
\}$ is orthonormalized base of subspace ${\sf m}\subset{\sf g}$ relative to
metric $(\cdot,\cdot)$ under relevant identifications. Evidently, in the base
$$\begin{array}{l}(\o g_{ij})=\left(\begin{array}{rr}
(1+\lambda)I_n&\mu I_n\\\mu I_n&(1-\lambda)I_n\end{array}\right);\\
(\o g^{ij})=\frac1{1-\l^2-\mu^2}\left(\begin{array}{rr}
(1-\lambda)I_n&\-\mu I_n\\-\mu I_n&(1+\lambda)I_n\end{array}\right);
\end{array}$$
where $I_n$ is the unit matrix of order $n$. Therefore, we have:
\begin{eqnarray*}
r(X)\eq\frac1{1-\l^2-\mu^2}\sum_{a=1}^n\,\{\,(1-\l)R(X,e_a)e_a+(1+\l)R(X,
e_{\hat a})e_{\hat a}-\\
&{}&-\mu R(X,e_a)e_{\hat a}-\mu R(X,e_{\hat a})e_a\,\}.
\end{eqnarray*}
A direct calculation in view of (36) shows that
\begin{eqnarray}
r(X)\eq\frac1{1-\l^2-\mu^2}\sum_{a=1}^n\,\{\,(-\frac{1-\l}4
+\frac{\mu^2(\mu-\l+1)}{4(1-\l^2-\mu^2)})[[X,e_a],e_a]+\nonumber\\
&{}&+(-\frac{1+\l}4+\frac{\mu^2(\l+\mu+1)}{4(1-\l^2-\mu^2)})[[X,
e_{\hat a}],e_{\hat a}]+\nonumber\\
&{}&+\frac{\mu(-2\mu^2-\l^2+3\l\mu-3\mu+2\l-1)}{4(1-\l^2-\mu^2)}
[[JX,e_a],e_a]-\\
&{}&-\frac{\mu(\l^2+2\mu^2+3\l\mu+2\l+3\mu+1)}{4(1-\l^2-\mu^2)}
[[JX,e_{\hat a}],e_{\hat a}].\nonumber
\end{eqnarray}
\procl{Lemma 2} Let $G$ is semisimple Lie group. Then the equality is valid:
$$g^{ij}[[X,e_i],e_j]=-X,$$
where $\{\,e_1,\ldots,e_n\,\}$ is some base of its Lie algebra $\sf g$,
$(g^{ij})$ is a contravariant metric tensor for Killing metric $\cal K$.
\eprocl
\demo{Proof} We have by definition: ${\cal K}(X,Y)=-\tr(\ad X\circ\ad Y),
\quad X,Y\in{\sf g}$.
Evidently, if $L:V\rightarrow V$ is any endomorphism of Euclidean space
$(V,\la\cdot,\cdot\ra)$, then $\tr L=g^{ij}\la L(e_i),e_j\ra$. Indeed,
$g^{ij}\la L(e_i),e_j\ra=g^{ij}\,L^k{}_i\,\la e_k,e_j\ra=g^{ij}\,g_{kj}\,
L^k{}_i=L^i{}_i$. Therefore,
\begin{eqnarray*}
{\cal K}(X,Y)\eq-g^{ij}\la[X,[Y,e_i]],e_j\ra=g^{ij}\la[Y,e_i],[X,e_j]\ra=\\
\eq-g^{ij}\la Y,[[X,e_j],e_i]\ra=-\la g^{ij}[[X,e_i],e_j],Y\ra\end{eqnarray*}
and hence
${\cal K}(X)=-g^{ij}[[X,e_i],e_j]$ is the result of subscript raising of Killing
form by itself, i.e. ${\cal K}=\id$ and hence $g^{ij}[[X,e_i],e_j]=-X$.\qed
\edemo
\procl{Corollary} Let $\{\,e_i\,\}$ be an orthonormalized (with respect to
Killing form) base of semisimple Lie algebra. Then
$$\sum_i[[X,e_i],e_i]=-X.\qeds$$
\eprocl
We return to considering metric $\la\cdot,\cdot\ra$ on manifold $M\times M$,
where $M=G$ is semisimple Lie group equipped with  Killing metric. Introduce
the notations
\begin{eqnarray*}
A\eq-\frac{1-\l}4+\frac{\mu^2(\mu-\l+1)}{4(1-\l^2-\mu^2)};\\
B\eq-\frac{1+\l}4+\frac{\mu^2(\mu+\l+1)}{4(1-\l^2-\mu^2)};\\
C\eq\frac{\mu(-2\mu^2-\l^2+3\l\mu-3\mu+2\l-1)}{4(1-\l^2-\mu^2)};\\
D\eq-\frac{\mu(2\mu^2+\l^2+3\l\mu+3\mu+2\l+1)}{4(1-\l^2-\mu^2)}.
\end{eqnarray*}
Let $\la\cdot,\cdot\ra$ be Einsteinian metric, i.e. $r=\e\id$. According to
(37) and Lemma 2 we have in this case: $A=B$, $C=D=0$ and hence $C-D=0$, i.e.
$\l\mu(3\mu+2)=0$. We the have following possibilities:
\begin{enumerate}
\item $\mu=0$. Then $C=D=0$, $A=-\frac{1-\l}4$, $B=-\frac{1+\l}4$. Therefore
$\l=0$, i.e. $\la X,Y\ra={\cal K}(X,Y)$.
\item $\mu\ne 0$, $\l=0$. Then $2\mu^2+3\mu+1=0$, i.e. $(2\mu+1)(\mu+1)=0$ and
hence $\mu=-\frac12$. Therefore in this case $\l=0$, $\mu=-\frac12$, i.e.
$\la X,Y\ra={\cal K}(X,Y)-\frac12(X,JY)$.
\item $\mu\ne 0$, $\l\ne 0$. Then $3\mu+2=0$, i.e. $\mu=-\frac23$, $\l^2=\frac19$;
$\l=\pm\frac13$. So, in this case $\l=\pm\frac13$, $\mu=-\frac23$, i.e.
$\la X,Y\ra={\cal K}(X,Y)+\frac13(X,IY)-\frac23(X,JY)$ or
$\la X,Y\ra={\cal K}(X,Y)-\frac13(X,IY)-\frac23(X,JY)$.
\end{enumerate}
\par
Inversely, if $\l=\mu=0$, then $A=B=-\frac14$, $C=D=0$ and hence $r=\frac14\id$.
If $\l=0$, $\mu=-\frac12$, then $A=B=-\frac5{24}$, $C=D=0$, $r=\frac5{18}\id$.
Finally, if $\l=\pm\frac13$, $\mu=-\frac23$, then $A=B=-\frac16$, $C=D=0$,
$r=\frac38\id$.
\par
So, we proved the following theorem giving a new method of building concrete
Einsteinian metrics:
\procl{Theorem 22} Let $G$ be semisimple Lie group equipped with Killing
metric $\cal K$. Then canonical almost antiquaternionic structure on manifold
$G\times G$ defines a two-parameter family of pseudo-Riemannian metrics
$$
\la X,Y\ra={\cal K}(X,Y)+\l{\cal K}(X,IY)+\mu{\cal K}(X,JY);\quad\l,\mu\in
{\bf R};\quad\l^2+\mu^2\ne 1.
$$
The metrics are Einsteinian in the following four cases only:
$$\begin{array}{lclrcl}
1)\ \l\eq\mu=0;&\quad\la X,Y\ra\eq{\cal K}(X,Y);\\
2)\ \l\eq0,\mu=-\frac12;&\quad\la X,Y\ra\eq{\cal K}(X,Y)-\frac12{\cal K}(X,JY);\\
3)\ \l\eq\frac13,\mu=-\frac23;&\quad\la X,Y\ra\eq{\cal K}(X,Y)+\frac13{\cal K}
(X,IY)-\frac23{\cal K}(X,JY);\\
4)\ \l\eq-\frac13,\mu=-\frac23;&\quad\la X,Y\ra\eq{\cal K}(X,Y)-\frac13
{\cal K}(X,IY)-\frac23{\cal K}(X,JY).\qeds
\end{array}$$
\eprocl
\procl{Corollary} Proper nearly Kaehlerian manifold $(G\times G,{\cal J}_0,
\la\cdot,\cdot\ra_0)$ that is built according to Theorem 21 for semisimple Lie
group $M=G$ by Killing metric is Einsteinian manifold.\qed
\eprocl

\end{document}